\documentclass[11pt]{amsart}
\input xy
\xyoption{all}
\usepackage{amsmath,amsthm, amscd, amssymb, amsfonts}

\usepackage{tikz} 
\usepackage{tikz-cd}

\newcommand{\eq}{\normalcolor{}}

\newcommand{\Hc}{{\mathcal H}}

\newcommand{\coev}{\mbox{coev}}
\newcommand{\ev}{\mbox{ev}}

\newcommand{\otk}{{\otimes_{\ku}}}
\newcommand{\botc}{{\boxtimes_\ca}}

\newcommand{\Mo}{{\mathcal M}}
\newcommand{\No}{{\mathcal N}}

\newcommand{\moda}{{\mathfrak m}}

\newcommand{\ca}{{\mathcal C}}

\newcommand{\ot}{{\otimes}}

\newcommand{\op}{\rm{op}}

\newcommand{\Ec}{{\mathcal E}}
\newcommand{\Ac}{{\mathcal A}}

\newcommand{\ere}{{\mathcal R}}
\newcommand{\ele}{{\mathcal L}}

\newcommand{\Do}{{\mathcal D}}

\newcommand{\Bc}{{\mathcal B}}

\newcommand{\Ker}{\mbox{\rm Ker\,}}
\newcommand{\cok}{\mbox{\rm coKer\,}}
\newcommand{\rev}{\rm{rev}}

\newcommand{\bop}{\rm{bop}}

\newcommand{\btc}{{\boxtimes_{\C}}}

\newcommand{\ku}{{\Bbbk}}

\newcommand{\uno}{ \mathbf{1}}
\newcommand{\C}{{\mathcal C}}
\newcommand{\sy}{{\mathbb S}}

\newcommand{\id}{\mbox{\rm id\,}}

\newcommand{\Id}{\mbox{\rm Id\,}}

\newcommand{\Res}{\mbox{\rm Res\,}}

\newcommand{\vect}{\mbox{\rm vect\,}}

\newcommand{\Nat}{\mbox{\rm Nat\,}}
\newcommand{\Rex}{\mbox{\rm Rex\,}}
\newcommand{\Fun}{\operatorname{Fun}}

\newcommand\bal{\operatorname{Bal}}

\newcommand\Hom{\operatorname{Hom}}
\newcommand\uhom{\underline{\Hom}}

\newcommand{\End}{\operatorname{End}}

\theoremstyle{plain}

\numberwithin{equation}{section}

\newtheorem{teo}{Theorem}[section]

\newtheorem{lema}[teo]{Lemma}

\newtheorem{cor}[teo]{Corollary}

\newtheorem{prop}[teo]{Proposition}

\newtheorem{claim}{Claim}[section]

\theoremstyle{definition}

\newtheorem{defi}[teo]{Definition}

\theoremstyle{remark}

\newtheorem{rmk}[teo]{Remark}

\def\pf{\begin{proof}}

\def\epf{\end{proof}}

\theoremstyle{remark}

\begin{document}

\title[  (co)ends  for representations of tensor categories]
{ (co)ends  for representations of tensor categories}
\author[    Bortolussi and Mombelli  ]{ Noelia Bortolussi and Mart\'in Mombelli
 }

\address{Facultad de Matem\'atica, Astronom\'\i a y F\'\i sica
\newline \indent
Universidad Nacional de C\'ordoba
\newline
\indent CIEM -- CONICET
\newline \indent Medina Allende s/n
\newline
\indent (5000) Ciudad Universitaria, C\'ordoba, Argentina}
 \email{ bortolussinb@gmail.com, nbortolussi@famaf.unc.edu.ar }
\email{martin10090@gmail.com, mombelli@mate.uncor.edu
\newline \indent\emph{URL:}\/ https://www.famaf.unc.edu.ar/$\sim$mombelli}

\begin{abstract} We generalize the notion of ends and coends in category theory to  the realm of module categories over finite tensor categories. We call this new concept \textit{module (co)end}. This tool allows us to give different proofs to several known results in the theory of representations of finite tensor categories. As a new application, we present a description of the relative Serre functor for module categories in terms of a module coend, in a analogous way as a Morita invariant description of the Nakayama functor of abelian categories presented in \cite{FSS0}.
\end{abstract}

\date{\today}
\maketitle

%\tableofcontents

\section*{Introduction}
Throughout this paper, $\ku$ will denote a field,  all categories will be finite (in the sense of \cite{EO}) abelian $\ku$-linear categories, and all functors will be additive $\ku$-linear. Given categories $\Mo, \Ac$, and a functor $S:\Mo^{\op}\times \Mo\to \Ac$ the notion of the end $\int_{M\in \Mo} S$ and coend $\int^{M\in \Mo} S$ is a standard and very useful concept in category theory. The end of the functor $S$ is an object  $\int_{M\in \Mo} S\in\Ac$  together  with dinatural transformations $$\pi_M:\int_{M\in \Mo} S \xrightarrow{\,..\,}  S(M,M)$$ with the following universal property; for any pair $(B, d)$ consisting of an object $B\in \Ac$ and a dinatural transformation $d_M: B \xrightarrow{..} S(M,M)$, there exists a \textit{unique} morphism $h: B\to \int_{M\in \Mo} S$ in $\Ac$ such that 
$$d_M = \pi_M \circ h \quad \text{ for any } M \in \Mo.$$

The notion of coend is defined dually.
\medbreak

If  $\Mo$ is a finite abelian, $\ku$-linear category, $\Mo$ can be thought of as a module category over  $\vect_\ku$, the tensor category of finite dimensional vector $\ku$-spaces. If $\Mo=\moda_A$ is the category of finite dimensional right $A$-modules, where $A$ is a finite dimensional $\ku$-algebra, then $\Mo$ has a left $\vect_\ku$-action 
 $$\vect_\ku\times \moda_A\to  \moda_A$$
 $$(V, M)\mapsto V\otk M, $$
where the right action on $V\otk M$ is given on the second tensorand. If $S: (\moda_A)^{\op}\times \moda_A\to \Ac$ is any  functor, it posses a canonical natural isomorphism
$$\beta^V_{M,N}: S(M, V\otk N)\to S(V^*\otk M,N),$$
for any $V\in \vect_\ku$, $M,N\in \moda_A$. See Proposition \ref{restriction-vect}. The existence of $\beta$ essentially follows from the additivity of the functor $S$.

 If in addition $p_M:E\xrightarrow{..} S(M,M)$ is any dinatural transformation, it satisfies equation
\begin{equation}\label{dwq1} S(\ev_V\otk \id_M,  \id_M) p_M= S(m_{V^*,V,M},  \id_M) \beta^V_{V\otk M,M} p_{V\otk M},
\end{equation}
for any $V\in \vect_\ku$. This is proven also in Proposition \ref{restriction-vect}. This equation follows from the dinaturality of $p$. Here $\ev_V:V^*\otk V\to \ku$ is the evaluation map, and $m_{W,V,M}:(W\otk V)\otk  M\to W\otk (V \otk  M)$ is the canonical associativity of vector spaces. This implies that the end of $S$ is the universal object among all dinatural transformations that satisfy \eqref{dwq1}. A similar observation can be made for the coend. This is the starting point to generalize the notion of (co)end, where we will replace the category $\vect_\ku$ with an arbitrary tensor category.

\medbreak

Let $\ca$ be a  tensor category, and $\Mo$ be a left $\ca$-module category with action given by  $\rhd: \ca \times \Mo \to \Mo$.  This action induces a right action of $\ca$ on the opposite category $\Mo^{\op}$:
$$ \blacktriangleleft :\Mo^{\op}\times \ca\to \Mo^{\op}.$$
\begin{equation*}
\overline{M}\blacktriangleleft X=\overline{X^*\triangleright M},
\end{equation*}
Here $\overline{M}$ is the object $M$ thought as an object in $\Mo^{\op}.$ Assume $S:\Mo^{\op}\times \Mo\to \Ac$ is a functor. We can  produce then two functors:
$$S\circ (\Id\times \triangleright),\,\, S\circ (\blacktriangleleft\times \Id): \Mo^{\op}\times \ca\times  \Mo\to \Ac.$$
Assume there exists a natural isomorphism $\beta: S\circ (\Id\times \triangleright)\to S\circ (\blacktriangleleft\times \Id)$, that is
 $$\beta^X_{M,N}: S(M,X\triangleright N)\to S(X^*\triangleright M,N).$$
We call this isomorphism a \textit{pre-balancing} of $S$. In this general case, the pre-balancing is an extra structure of the functor $S$. We define the \textit{module end} of $S$
to be an object  $E \in \Ac$
 that comes with dinatural transformations $\pi_M: E \xrightarrow{ . .} S(M,M)$ such that the  equation
\begin{equation}\label{dwq2} S(\ev_X\triangleright \id_M,  \id_M) \pi_M= S(m_{X^*,X,M},  \id_M) \beta^X_{X\triangleright M,M} \pi_{X\triangleright M},\end{equation}
is fulfilled, and it is universal among all objects in $\Ac$ with dinatural transformations that satisfy \eqref{dwq2}. Unlike the case $\ca= \vect_\ku$, it may happen that a dinatural transformation does not satisfy  \eqref{dwq2}.
We denote the module end as $\oint_{M\in \Mo} (S, \beta),$ or sometimes simply as  $\oint_{M\in \Mo} S$ whenever the pre-balancing $\beta$ is undertstood from the context.

An analogous definition can be made to define \textit{module coend}, and also to define  module ends and coends starting from \textit{right} $\ca$-module categories.
\medbreak

In Section \ref{Section:mcoends} we introduce the module (co)ends, and we prove several results that extend known properties of  (co)ends. We prove that, when the tensor category $\ca=\vect_\ku$ our definition coincides with the usual (co)ends. See Proposition \ref{restriction-vect}. We also study what happens when we restrict the module (co)ends to a tensor subcategory. See Proposition \ref{restriction-sucategory}.

In Section \ref{Section:application} we give several applications. If $\Mo, \No$ are left $\ca$-module categories, and $F, G:\Mo\to \No$ are $\ca$-module functors, the functor 
$$ \Hom_\No(F(-),G(-)): \Mo^{\op}\times \Mo\to \vect_\ku$$
has a canonical pre-balancing $\gamma$, and we prove that, there is an isomorphism 
$$\Nat_{\!m}(F,G)\simeq \oint_{M\in \Mo}  (\Hom_\No(F(-),G(-)), \gamma).$$
Here $\Nat_{\!m}(F,G)$ is the space of natural \textit{module} transformations between $F$ and $G$. See Proposition \ref{end-natural-transf}. Using this result we can set up a triangle of adjoint equivalences of categories
 $$ \begin{tikzcd}[row sep=11ex]
  ~ & \Mo^{\op} \btc \No ~ \ar{dl}[xshift=-2pt]{L_{\Mo,\No}} \ar[xshift=-2pt]{dr}[swap]{\widetilde{L}_{\Mo,\No} }
  & ~ \\
  \Fun_\ca(\Mo^{\bop} ,\No)\ar[yshift=3pt]{rr}[yshift=2pt]{\Theta_{\Mo,\No}} \ar[xshift=-12pt]{ur}[xshift=2pt]{\chi_{\Mo,\No}}
    & ~ & \Fun_\ca(\Mo ,\No)\ar[xshift=10pt]{ul}[swap]{ \Upsilon_{\Mo,\No}},
  \end{tikzcd}
 $$
generalizing the triangle presented in \cite{FSS0}. Here it is required that $\Mo, \No$ are \textit{exact} module categories. Here $\Mo^{\op}$ is the opposite category endowed with a \emph{right} $\ca$-action that comes from the action of $\Mo$ twisted by a (right) dual. Also $\Mo^{\bop}=(\Mo^{\op})^{\op}.$ See Subsection \ref{bimodule-categories} for these definitions. Observe that, $\Mo^{\bop}=\Mo$ as categories, but as a $\ca$-module category $\Mo^{\bop}$ has the deformed action of $\Mo$ by a double dual.

The equivalences presented above are:
$$L_{\Mo,\No}: \Mo^{\op}\boxtimes_\ca \No \to \Fun_\ca(\Mo^{\bop}, \No),$$
$$M\btc N \mapsto \uhom_{\Mo^{\op}}( -, M)\triangleright N,$$
$$ \chi_{\Mo,\No}: \Fun_\ca(\Mo^{\bop}, \No)\to \Mo^{\op}\boxtimes_\ca \No,$$
$$F\mapsto   \oint^{\overline{U}\in \Mo^{\op}} \overline{U}\btc F(\overline{\overline{U}}),$$
and on the other side of the triangle we have equivalences
$$\widetilde{L}_{\Mo,\No}: \Mo^{\op}\boxtimes_\ca \No \to \Fun_\ca(\Mo, \No)$$
$$M\btc N \mapsto \uhom_{\Mo^{\op}}( M, -)^*\triangleright N,$$
$$ \Upsilon_{\Mo,\No}: \Fun_\ca(\Mo ,\No)\to \Mo^{\op}\boxtimes_\ca \No,$$ 
$$ F\mapsto   \oint_{M\in \Mo}  \overline{M}\boxtimes_\ca F(M).$$
Here $\Theta_{\Mo,\No}=\widetilde{L}_{\Mo,\No}\circ \chi_{\Mo,\No}. $ 
The proof that $ \widetilde{L}$ and $ L$ are equivalences is given in Lemma \ref{aboutL}. The explicit description of a quasi-inverse of $ \widetilde{L}$ is given in Theorem \ref{inv-L}.

As a consequence of these equivalences, in Corollary \ref{pw-bim} we obtain a kind of Peter-Weyl theorem for the regular $A$-bimodule $A={}_A A_A$; that is, if $A\in \ca$ is an algebra such that, the module category $\ca_A$ is exact, there is an isomorphism of $A$-bimodules:
$$A\simeq \oint_{M\in \Mo} {}^*M \otimes M. $$
We also prove that the module functor $\Theta_{\Mo,\Mo^{\bop}}(\Id)$ is equivalent, as module functors, to the (right) relative Serre functor of $\Mo$. See Theorem \ref{formula-serre-relat}. This description is an analogous form of the Morita invariant description of the Nakayama functor presented in \cite{FSS0}.

If $\ca$ and $\Do$ are Morita-equivalent tensor categories, this means that there exists an invertible $(\ca, \Do)$-bimodule category $\Bc$; we prove that the correspondence
$$ \Mo \mapsto \Fun_\ca(\Bc, \Mo), \quad \No \mapsto \Fun_{\Do}(\Bc^{\op}, \No)$$
is in fact part of a 2-equivalence between the 2-categories of $\ca$-module categories and $\Do$-module categories. This result was proven in \cite{EO}. We show in Theorem \ref{bij-corresp} that, for any $\Do$-module category $\No$, the functor
$$ \Fun_{\ca}(\Bc, \Fun_{\Do}(\Bc^{\op}, \No)) \to \No$$
$$H\mapsto \oint_{B\in \Bc} H(B)(B)$$
is an equivalence of $\Do$-module categories.
\medbreak

In the last Section we show that, the functor  $\Upsilon:(\ca^*_\Mo)^*_\Mo \to \ca$ defined as 
$$\Upsilon(G)= \oint_{M\in \Mo}  \uhom(M, G(M))$$ 
is a quasi-inverse of the canonical functor $$can: \ca\to  (\ca^*_\Mo)^*_\Mo,
\quad can(X)(M)=X\rhd M. $$ See Theorem \ref{ddual-tensor}.
\subsection*{Acknowledgments} This  work  was  partially supported by CONICET and Secyt (UNC),  Argentina. We would like to thank the referee for his/her many comments that significantly improved  the presentation of the paper. In particular, Proposition \ref{ach} was suggested by him/her.

\subsection*{Preliminaries and Notation}

We  denote by  $\vect_\ku$ the category of finite-dimensional $\ku$-vector spaces.  If $\Mo, \No$ are categories, and $F:\Mo\to \No$ is a functor, we shall denote by $F^{r.a.}, F^{l.a.}:\No\to \Mo$ its right and left adjoint, respectively.

For any category $\Mo$, the opposite category will be denoted   by $\Mo^{\op}$.  We shall denote by $\overline{M},  \overline{f}$ objects and morphisms in $\Mo^{\op}$ that correspond to $M$ and $f$. We shall also denote by $F^{\op}:\Mo^{\op}\to \No^{\op}$ the opposite functor to $F$; that is, the functor defined as $F^{\op}(\overline{M})=\overline{F(M)} $, $F^{\op}(\overline{f})=\overline{F(f)} $ for any object $M$ and any morphism $f$.

 \section{Finite tensor categories}\label{Section:tensorcat}
  For basic notions on finite tensor categories we refer to \cite{EGNO}, \cite{EO}. Let $\C$ be a  tensor category over $\ku$; that is a rigid monoidal category with simple unit object $\uno$.
  
    If $\ca$ has  associativity constraint given by $$a_{X,Y,Z}:(X\ot Y)\ot Z\to X\ot (Y\ot Z), $$ we shall denote by $\ca^{\rev}$, the tensor category whose underlying abelian category is $\ca$, with \textit{reverse} monoidal product 
    $$\ot^{\rev}:\ca\times \ca\to \ca, \quad X\ot^{\rev }Y=Y\ot X,$$ and associativity constraints
 $$a^{\rev}_{X,Y,Z}:( X\ot^{\rev }Y) \ot^{\rev} Z\to  X\ot^{\rev }(Y \ot^{\rev} Z),$$
  $$a^{\rev}_{X,Y,Z}:= a^{-1}_{Z,Y,X},$$
  for any $X,Y, Z\in \ca$.
  It is well known that for any pair of objects $X, Y\in \ca$ there are canonical isomorphisms
\begin{equation}\label{duals-tensor-product}
\begin{split} \phi^r_{X,Y}: (X\ot Y)^* \to Y^* \ot X^*,\\
\phi^l_{X,Y}:  {}^*(X\ot Y)\to {}^*Y\ot  {}^*X.
\end{split}
\end{equation} 

For any $X\in \ca$ we shall denote by 
$$ \ev_X:X^*\ot X\to \uno, \quad \coev_X:\uno\to X\ot X^*$$
the evaluation and coevaluation. Abusing of the notation, we shall also denote by
$$\ev_X:X\ot {}^*X\to \uno,  \quad \coev_X:\uno\to  {}^*X\ot X$$
the evaluation and coevaluation for the left duals.
If $f:X\to Y$ is an isomorphism in $\ca$ then
\begin{align}\label{duality-morph} \ev_Y (f\ot  \id_Y)=\ev_X (\id_X\ot {}^*f).
\end{align}
For any $X, Y\in \ca$ the following identities hold
\begin{equation}\label{ev-tensor-prod}\begin{split}
\ev_{X\ot Y}=\ev_X (\id_X\ot \ev_Y\ot \id_{{}^*X})(\id_{X\ot Y}\ot \phi^l_{X,Y}),\\
( \phi^l_{X,Y} \ot \id_{X\ot Y}) \coev_{X\ot Y}= (\id_{{}^*Y}\ot \coev_X\ot \id_Y) \coev_Y.
\end{split}
\end{equation}
Off course that similar identities hold for the right duals, but they won't be needed.
  \subsection{Algebras in tensor categories}
  
 In this subsection we assume that $\ca$ is a strict tensor category, this means in particular that the associativity constraints are the identities.  Let  $A, B\in \ca$ be  algebras. We shall denote by $$\ca_A, \, \,   {}_A\ca,  \,   \,  {}_A\ca_B$$ the categories of right $A$-modules, left $A$-modules and $(A, B)$-bimodules in $\ca$, respectively.  If $V\in \ca_A$ is a right $A$-module with action given by $\rho_V:V\ot A\to V$, and  $W\in{}_A\ca$ is a left $A$-module with action given by $\lambda_W:A\ot W\to W$,  we shall denote by $\pi^A_{V,W}: V \ot W\to V\ot_A W$  the coequalizer of the maps
$$\rho_V\ot\id_W,\,  \id_V\ot \lambda_W:(V \ot A)\ot W\longrightarrow V\ot W.$$
 An object in the category ${}_A\ca_B$ will be denoted as $(V, \lambda_V, \rho_V) \in {}_A \ca_B$, where  $ \lambda_V:A\ot V\to V$ is the left action, and $\rho_V:V\ot B\to V$ is the right action. Since the tensor product is exact in both variables, then 
 $$ \pi^A_{V, W\ot U}=\pi^A_{V,W}\ot \id_U,$$
 for any  $V\in \ca_A, W\in  {}_A\ca$, $U\in \ca$. We are going to freely use this fact without further mention.
 
  \begin{lema}\label{dual-module1} Assume that $\ca$ is a  tensor category and $A, B\in \ca$ are algebras. The following statements hold:
  \begin{itemize}
  \item[(i)] If $M\in \ca_A$ then ${}^*M\in {}_A\ca$.
  
   \item[(ii)] There are natural isomorphisms 
   \begin{equation}\label{hom-iso-1}
   \Hom_{B}(M\ot_A V, U)\simeq \Hom_{(A,B)}(V, {}^*M\ot U),
  \end{equation} 
  \begin{equation}\label{hom-iso-2}
   \Hom_{A}(M, X\ot N)\simeq \Hom_\ca(M\ot_A {}^*N, X),
  \end{equation} 
  \begin{equation}\label{hom-iso-3}
   \Hom_{A}(M, X\ot N)\simeq  \Hom_{A}(X^*\ot M, N),
  \end{equation} 
      for any $X\in \ca, M, N\in \ca_A, V\in {}_A\ca_B, U\in \ca_B. $
  \end{itemize}
  \end{lema}
  \pf (i). If $M\in \ca_A$ then  ${}^*M $ has structure of left $A$-module via $\lambda_{{}~M} : A \ot {}^*M \rightarrow {}^*M$ defined as \begin{equation}\label{dual-leftA}
  \lambda_{{}^*M} = (\id_{{}^*M} \ot \ev_M)(\id_{{}^*M} \ot \rho_M \ot \id_{{}^*M})(\coev_M \ot \id_{A \ot {}^*M}).
  \end{equation}
  
  (ii). Let us prove only the first isomorphism. The others follow similarly. The object $M\ot_A V$ has a right $B$-module structure as follows. Consider $\phi: M \ot V \ot B \rightarrow M \ot_A V$, $\phi = \pi_{M,V} (\id_M \ot \rho_V)$. Then, $\rho_{M\ot_A V}: M \ot_A V \ot B \rightarrow M \ot_A V$ is defined as the unique morphism such that
  \begin{equation}\label{def-rho}
      \rho_{M \ot_A V} (\pi_{M,V} \ot \id_B) = \phi.
  \end{equation}
  Define $\Phi: \Hom_B (M \ot_A V, U) \rightarrow \Hom_{(A,B)}(V, {}^*M \ot U)$ as
 \begin{equation}\label{def-phi-b}
 \Phi(f)=(\id_{{}^*M} \ot f \pi_{M,V})(\coev_M \ot \id_V),
 \end{equation}
 for any $f\in \Hom_B (M \ot_A V, U).$  Let us show $\Phi(f)$ is a morphism of  $(A,B)$-bimodules. We need to prove that
  \begin{equation}\label{A-mod eq}
      (\lambda_{{}^*M} \ot \id_U)(\id_A \ot \Phi(f)) = \Phi(f) \lambda_V,
  \end{equation}
  and
  \begin{equation}\label{B-mod eq}
      (\id_{{}^*M} \ot \rho_U)(\Phi(f) \ot \id_B) = \Phi(f) \rho_V,
  \end{equation}
   for any $(M,\rho_M) \in \ca_A$, $(V, \lambda_V, \rho_V) \in {}_A \ca_B$ and $(U, \rho_U) \in \ca_B$. Here $   \lambda_{{}^*M} $ is the left action of $A$ on ${}^*M$ presented in \eqref{dual-leftA}.
   
The left hand side of \eqref{A-mod eq} is equal to 
   \begin{align*}
       (\lambda_{{}^*M} & \ot \id_U)(\id_A \ot \Phi(f)) =\\
       & = (\id_{{}^*M} \ot \ev_M \ot \id_U)(\id_{{}^*M} \ot \rho_M \ot \id_{{}^*M\ot U})(\coev_M \ot \id_{A \ot {}^*M \ot U})\\
       & (\id_{A \ot {}^*M} \ot f \pi_{M,V})(\id_A \ot \coev_M \ot \id_V)\\
       & = (\id_{{}^*M} \ot \ev_M \ot \id_U)(\id_{{}^*M} \ot \rho_M \ot \id_{{}^*M \ot U})(\id_{{}^*M \ot M \ot A \ot {}^*M} \ot f \pi_{M,V})\\
       & (\coev_M \ot \id_{A \ot {}^*M \ot M \ot V})(\id_A \ot \coev_M \ot \id_V)\\
       & = (\id_{{}^*M} \ot \ev_M \ot \id_V)(\id_{{}^*M \ot M \ot {}^*M} \ot f \pi_{M,V})(\id_{{}^*M} \ot \rho_M \ot \id_{{}^*M \ot M \ot V})\\
       & (\id_{{}^*M \ot M \ot A} \ot \coev_M \ot \id_V)(\coev_M \ot \id_{A \ot V})\\
       & = (\id_{{}^*M} \ot f \pi_{M,V})(\id_{{}^*M} \ot \ev_M \ot \id_{M \ot V})(\id_{{}^*M \ot M} \ot \coev_M \ot \id_V)\\
       & (\id_{{}^*M} \ot \rho_M \ot \id_V)(\coev_M \ot \id_{A \ot V})\\
       & =  (\id_{{}^*M} \ot f )(\id_{{}^*M} \ot \pi_{M,V}(\rho_M \ot \id_V))(\coev_M \ot \id_{A \ot V})\\
       & = (\id_{{}^*M} \ot f\pi_{M,V} )(\id_{{}^*M \ot M} \ot \lambda_V)(\coev_M \ot \id_{A \ot V})\\
       & = (\id_{{}^*M} \ot f\pi_{M,V} )(\coev_M \ot \id_V) \lambda_V\\
       & = \Phi(f) \lambda_V.
  \end{align*}
  The first equality is by the definition of $\lambda_{{}^*M}$ and $\Phi(f)$. The fifth equality follows from the rigidity axioms. The sixth equality is consequence of $\pi_{M,V}$ being the coequalizer  of  $\rho_M \ot \id_V, \id_M \ot \lambda_V$. The last equality follows by the definition of $\Phi(f)$.

 Since $f$ is a $B$-module morphism, 
  \begin{equation}\label{f B-mod}
      \rho_U (f \ot \id_B) = f \rho_{M \ot_A V}.
  \end{equation}
  Using \eqref{def-rho}, this equation implies 
  \begin{equation}\label{relation-f-rho}
      \rho_U (f \pi_{M,V} \ot \id_B) = f \rho_{M \ot_A V} (\pi_{M,V} \ot \id_B) = f \pi_{M,V} (\id_M \ot \rho_V).
  \end{equation}
  Let us prove \eqref{B-mod eq}. The left hand side of \eqref{B-mod eq} is equal to
  \begin{align*}
      (\id_{{}^*M} \ot \rho_U)(\Phi(f) \ot \id_B) & = (\id_{{}^*M} \ot \rho_U (f \pi_{M,V} \ot \id_B))(\coev_M \ot \id_{V \ot B})\\
      & = (\id_{{}^*M} \ot f \pi_{M,V}(\id_M \ot \rho_V))(\coev_M \ot \id_{V \ot B})\\
      & = (\id_{{}^*M} \ot f \pi_{M,V})(\coev_M \ot \id_V) \rho_V\\
      & = \Phi(f) \rho_V.
  \end{align*}
  The first equality is by the definition of $\Phi(f)$.
  The second equality follows from  \eqref{relation-f-rho}.
   And the last equality  again follows from the definition of $\Phi(f)$.
   
  Now, let us show that $\Phi$ has an inverse. Let us define $$\Psi: \Hom_{(A,B)} (V, {}^*M \ot U) \rightarrow \Hom_B (M\ot_A V, U)$$ as follows. Let $g \in \Hom_{(A,B)} (V, {}^*M \ot U)$. Define $\Psi(g)=h$ where $h: M \ot_A V \rightarrow U$ is the unique morphism such that 
 \begin{equation}\label{def-h}
     h \pi_{M,V} = (\ev_M \ot \id_U)(\id_M \ot g).
 \end{equation}
 Let us show $\Psi(g)$ is a $B$-module morphism. That is
 $$\rho_U (h \ot \id_B) = h \rho_{M \ot_A V}.$$
 For this, it is enough to prove
  $$\rho_U (h \pi_{M,V} \ot \id_B) = h \rho_{M \ot_A V}(\pi_{M,V} \ot \id_B).$$
  Starting from the left hand side
  \begin{align*}
   \rho_U (h \pi_{M,V} \ot \id_B) & =  \rho_U (\ev_M \ot \id_{U \ot B})(\id_M \ot g \ot \id_B)\\
   & = (\ev_M \ot\id_U)(\id_{M \ot {}^*M} \ot \rho_U)(\id_M \ot g \ot \id_B)\\
   & = (\ev_M \ot \id_U)(\id_M \ot g \rho_V)\\
   & = h \pi_{M,V} (\id_M \ot \rho_V)\\
   & = h \rho_{M \ot_A V} (\pi_{M,V} \ot \id_B).
  \end{align*}
  The first equality is by \eqref{def-h}.
  The third equality is consequence of $g$ being a $B$-module morphism.
  The fourth equality follows from \eqref{def-h} and the last equality follows from \eqref{def-rho}.
   Let us show $\Phi$ and $\Psi$ are inverses  of each another. Let  be $f \in \Hom_B (M \ot_A V, U)$. We have
  $$\Psi \Phi (f) = \Psi((\id_{{}^*M} \ot f \pi_{M,V})(\coev_M \ot \id_V)) = h$$
  where
  \begin{align*}
      h \pi_{M,V} & = (\ev_M \ot \id_U)(\id_{{}^*M \ot M} \ot f \pi_{M,V})(\id_M \ot \coev_M \id_V)\\
      & = f \pi_{M,V} (\ev_M \ot \id_{M \ot V})(\id_M \ot \coev_M \ot \id_V)\\
      & = f \pi_{M,V}
  \end{align*}
  
  The first equality is the definition of $h$, and the last equality follows from the rigidity axioms. Therefore, $h = f$ and $\Psi \Phi (f) = f$. The proof of $\Phi \Psi= \Id$ follows similarly.
 \medbreak
 
We shall only sketch the proof of isomorphism \eqref{hom-iso-2}. Define
\begin{equation}\label{defini-isos2}\begin{split}
 \Phi^A_{M,X,N}: \Hom_A(M, X\ot N)\to \Hom_\ca(M\ot_A {}^*N, X),\\
  \Phi^A_{M,X,N}(\alpha)\pi^A_{M, {}^*N}=(\id_X\ot \ev_N)(\alpha\ot\id_{ {}^*N}),
 \end{split}
\end{equation}
and 
 \begin{equation}\label{defini-isos22}\begin{split}
  \Psi^A_{M,X,N}:\Hom_\ca(M\ot_A {}^*N, X)\to  \Hom_A(M, X\ot N),\\
  \Psi^A_{M,X,N}(\alpha)=(\alpha \pi^A_{M, {}^*N}\ot\id_N)(\id_M\ot \coev_N).
    \end{split}
\end{equation}
It follows by a direct calculation that $\Phi^A_{M,X,N}$ and $\Psi^A_{M,X,N}$ are well-defined and they are one the inverse of the other.
  \epf

 \section{Representations of tensor categories} A  left \emph{module} category over  
$\ca$ is a  category $\Mo$ together with a $\ku$-bilinear 
bifunctor $\rhd: \ca \times \Mo \to \Mo$, exact in each variable,  endowed with 
 natural associativity
and unit isomorphisms 
$$m_{X,Y,M}: (X\otimes Y)\triangleright   M \to X\triangleright  
(Y\triangleright M), \ \ \ell_M: \uno \triangleright  M\to M.$$ 
These isomorphisms are subject to the following conditions:
\begin{equation}\label{left-modulecat1} m_{X, Y, Z\triangleright M}\; m_{X\otimes Y, Z,
M}= (\id_{X}\triangleright m_{Y,Z, M})\;  m_{X, Y\otimes Z, M}(a_{X,Y,Z}\triangleright\id_M),
\end{equation}
\begin{equation}\label{left-modulecat2} (\id_{X}\triangleright \ell_M)m_{X,{\bf
1} ,M}= r_X \triangleright \id_M,
\end{equation} for any $X, Y, Z\in\C, M\in\Mo.$ Here $a$ is the associativity constraint of $\C$.
Sometimes we shall also say  that $\Mo$ is a $\ca$-\emph{module category} or a representation of $\ca$.

\medbreak

Let $\Mo$ and $\Mo'$ be a pair of $\C$-modules. A\emph{ module functor} is a pair $(F,c)$, where  $F:\Mo\to\Mo'$  is a functor equipped with natural isomorphisms
$$c_{X,M}: F(X\triangleright M)\to
X\triangleright F(M),$$ $X\in  \ca$, $M\in \Mo$,  such that
for any $X, Y\in
\ca$, $M\in \Mo$:
\begin{align}\label{modfunctor1}
(\id_X \triangleright  c_{Y,M})c_{X,Y\triangleright M}F(m_{X,Y,M}) &=
m_{X,Y,F(M)}\, c_{X\otimes Y,M}
\\\label{modfunctor2}
\ell_{F(M)} \,c_{\uno ,M} &=F(\ell_{M}).
\end{align}

There is a composition
of module functors: if $\Mo''$ is a $\C$-module category and
$(G,d): \Mo' \to \Mo''$ is another module functor then the
composition
\begin{equation}\label{modfunctor-comp}
(G\circ F, e): \Mo \to \Mo'', \qquad  e_{X,M} = d_{X,F(M)}\circ
G(c_{X,M}),
\end{equation} is
also a module functor.

\smallbreak  

A \textit{natural module transformation} between  module functors $(F,c)$ and $(G,d)$ is a 
 natural transformation $\theta: F \to G$ such
that
\begin{gather}
\label{modfunctor3} d_{X,M}\theta_{X\triangleright M} =
(\id_{X}\triangleright \theta_{M})c_{X,M},
\end{gather}
 for any $X\in \ca$, $M\in \Mo$. The vector space of natural module transformations will be denoted by $\Nat_{\!m}(F,G)$. Two module functors $F, G$ are \emph{equivalent} if there exists a natural module isomorphism
$\theta:F \to G$. We denote by $\Fun_{\ca}(\Mo, \Mo')$ the category whose
objects are module functors $(F, c)$ from $\Mo$ to $\Mo'$ and arrows module natural transformations. 

\medbreak
Two $\C$-modules $\Mo$ and $\Mo'$ are {\em equivalent} if there exist module functors $F:\Mo\to
\Mo'$, $G:\Mo'\to \Mo$, and natural module isomorphisms
$\Id_{\Mo'} \to F\circ G$, $\Id_{\Mo} \to G\circ F$.
\medbreak
A module is
{\em indecomposable} if it is not equivalent to a direct sum of
two non trivial modules. Recall from \cite{EO}, that  a
module $\Mo$ is \emph{exact} if   for any
projective object
$P\in \ca$ the object $P\triangleright M$ is projective in $\Mo$, for all
$M\in\Mo$. If $\Mo$ is an exact indecomposable module category over $\ca$, the dual category $\ca^*_\Mo=\End_\ca(\Mo)$ is a finite tensor category \cite{EO}. The tensor product is the composition of module functors.

\medbreak
 A \emph{right module category} over $\ca$
 is a finite  category $\Mo$ equipped with an exact
bifunctor $\triangleleft:  \Mo\times  \ca\to \Mo$ and natural   isomorphisms 
$$\widetilde{m}_{M, X,Y}: M\triangleleft (X\ot Y)\to (M\triangleleft X) \triangleleft Y, \quad r_M:M\triangleleft \uno\to M$$ such that
\begin{equation}\label{right-modulecat1} \widetilde{m}_{M\triangleleft X, Y ,Z }\; \widetilde{m}_{M,X ,Y\ot Z } (\id_M \triangleleft a_{X,Y,Z})=
(\widetilde{m}_{M,X , Y}\triangleleft \id_Z)\, \widetilde{m}_{M,X\ot Y ,Z },
\end{equation}
\begin{equation}\label{right-modulecat2} (r_M\triangleleft \id_X)  \widetilde{m}_{M,\uno, X}= \id_M\triangleleft l_X.
\end{equation}

If $\Mo,  \Mo'$ are right $\ca$-modules, a module functor from $\Mo$ to $  \Mo'$ is a pair $(T, d)$ where
$T:\Mo \to \Mo'$ is a  functor and $d_{M,X}:T(M\triangleleft X)\to T(M)\triangleleft X$ are natural  isomorphisms
such that for any $X, Y\in
\ca$, $M\in \Mo$:
\begin{align}\label{modfunctor11}
( d_{M,X}\ot \id_Y)d_{M\triangleleft X, Y}T(m_{M, X, Y}) &=
m_{T(M), X,Y}\, d_{M, X\ot Y},
\\\label{modfunctor21}
r_{T(M)} \,d_{ M,\uno} &=T(r_{M}).
\end{align}

The next result is well-known. See for example \cite[Corollary 2.13.]{DSS}, \cite[Prop. 2.2.4]{GSch2}.
\begin{lema}\label{adjoint-module-functor} Let $\Mo, \No$ be left $\ca$-module categories, and $F, G:\Mo\to \No$ are $\ca$-module functors. 
\begin{itemize}
\item[(i)] The right and left adjoint of $F$, if they exist, have structure of $\ca$-module functor. 

\item[(ii)] If $F\simeq G$ as $\ca$-module functors, then $F^{l.a}\simeq G^{l.a}$, $F^{r.a}\simeq G^{r.a}$ as $\ca$-module functors.

\item[(iii)] If $F_1, F_2$ are composable $\ca$-module functors, there exists an isomorphism of $\ca$-module functors 
$$(F_1\circ F_2)^{l.a}\simeq F^{l.a}_2\circ F^{l.a}_1, \,  (F_1\circ F_2)^{r.a}\simeq F^{r.a}_2\circ F^{r.a}_1.$$ \qed
\end{itemize}
\end{lema}

\subsection{Bimodule categories}\label{bimodule-categories}
Assume that $\ca, \Do, \Ec$ are   tensor categories. A $(\ca, \Do)-$\emph{bimodule category}  is a category $\Mo$  with left $\ca$-module category structure 
$\triangleright: \ca\times \Mo\to \Mo$,
 and right $\Do$-module category  structure $ \triangleleft: \Mo\times \Do \to \Mo$,
equipped with natural
isomorphisms $$\{\gamma_{X,M,Y}:(X\triangleright M) \triangleleft Y\to X\triangleright  (M\triangleleft Y), X\in\ca, Y\in\Do, M\in \Mo\}$$ satisfying 
certain axioms. For details the reader is referred to \cite{Gr}, \cite{Gr2}.

If $\Mo$ is a   right $\ca$-module category then the opposite category $\Mo^{\op}$ has a left $\ca$-action given by 
$$\ca\times\Mo^{\op} \to \Mo^{\op},$$ $$(X, M)\mapsto  M\triangleleft\, X^*,$$ and associativity
 isomorphisms $m^{\op}_{X,Y,M}= m_{M,Y^*,X^*}(\id_M\triangleleft \phi^r_{X,Y}).$ Analogously, if $\Mo$ is a left  $\ca$-module category then $\Mo^{\op}$ has structure of right $\ca$-module category, with action given by 
$$\Mo^{\op}\times \ca\to  \Mo^{\op},$$   $$(M,X)\mapsto X^*\triangleright M,$$ 
 with associativity constraints 
 $m^{\op}_{M,X,Y}=m_{Y^*, X^*, M} (\phi^r_{X,Y}\triangleright\id_M)$ for all $X, Y\in \ca, M\in \Mo$.  If  $\Mo$ is a $(\ca,\Do)$-bimodule category then $\Mo^{\op}$ is a 
$(\Do,\ca)$-bimodule category. 

\medbreak

If $\Mo$ is a left  $\ca$-module category, we shall denote by $\Mo^{\bop}=(\Mo^{\op})^{\op}$. That is, $\Mo^{\bop}=\Mo$ as categories, but the left action of $\ca$ on  $\Mo^{\bop}$ is
$$\blacktriangleright:\ca\times  \Mo^{\bop}\to \Mo^{\bop},$$
$$X\blacktriangleright M=X^{**} \triangleright M, $$
for any $X\in \ca$, $M\in \Mo$.

\begin{rmk}\label{op-right} There is no problem to define the  actions on the category $\Mo^{\op}$  using left duals instead of right duals. Our choice of using right duals is related to the choice of functors $L, \widetilde{L}$ presented later in \eqref{L-equivalence1}, \eqref{L-equivalence2}.
\end{rmk}

Assume that $\Mo$ is a $(\ca, \Do)$-bimodule category, and   $\No$ is a  $(\ca,\Ec)$-bimodule category. The category $\Fun_\ca(\Mo, \No)$ has a structure of  $(\Do,\Ec)$-bimodule category.
Let us briefly describe this structure. For more details, the reader is referred to \cite{Gr}. The left and right actions are given by 
$$\triangleright: \Do \times \Fun_\ca(\Mo, \No)\to \Fun_\ca(\Mo, \No),$$
$$ \triangleleft: \Fun_\ca(\Mo, \No)\times \Ec \to \Fun_\ca(\Mo, \No), $$
where 
\begin{equation}\label{actions-funct-categ}
(X \triangleright F)(M)=F(M \triangleleft X),\quad  (F  \triangleleft Y)(M)=F(M)\triangleleft Y, 
\end{equation} 
for any $X\in \Do$,  $Y\in \Ec$, $F\in \Fun_\ca(\Mo, \No)$ and $M\in \Mo$.

\subsection{The internal Hom}\label{subsection:internal hom} Let $\ca$ be a  tensor category and $\Mo$ be  a left $\C$-module category. For any pair of objects $M, N\in\Mo$, the \emph{internal Hom} is an object $\uhom(M,N)\in \C$ representing the left exact functor $$\Hom_{\Mo}(-\triangleright M,N):\ca^{\op}\to \vect_\ku.$$ This means that, there are natural isomorphisms, one the inverse of each other,
\begin{equation}\label{Hom-interno}\begin{split}\phi^X_{M,N}:\Hom_{\ca}(X,\uhom(M,N))\to \Hom_{\Mo}(X\triangleright M,N), \\
\psi^X_{M,N}:\Hom_{\Mo}(X\triangleright M,N)\to \Hom_{\ca}(X,\uhom(M,N)),
\end{split}
\end{equation}
 for all $M, N\in \Mo$, $X\in\ca$. 
Sometimes we shall denote the internal Hom of the module category $\Mo$ by $\uhom_\Mo$ to emphasize that it is related to this module category. Similarly, if $\No$ is a right $\C$-module category, for any pair $M, N\in \No$ the internal hom  is the object $\uhom(M,N)\in \C$ representing the left exact functor $$\Hom_{\Mo}(M\triangleleft 
-,N):\ca^{\op}\to \vect_\ku.$$

\begin{lema}\label{internal-hom-actions} The following statements hold.
\begin{itemize}
\item[1.] Let $\Mo$ be a left $\ca$-module category. There are natural isomorphisms
$$ \uhom_{\Mo}(X\triangleright M, N)\simeq \uhom_{\Mo}( M, N)\ot X^*, $$
$$\uhom_{\Mo}( M, X\triangleright N)\simeq X\ot \uhom_{\Mo}( M, N).$$
for any $M, N\in \Mo$, $X\in \ca$.
\item[2.] Analogously, if $\No$ is a right $\ca$-module category, there are natural isomorphisms
$$ \uhom_{\No}(M \triangleleft X, N)\simeq  {}^*X\ot \uhom_{\No}( M, N), $$
$$\uhom_{\No}( M,  N\triangleleft X)\simeq \uhom_{\No}( M, N)\ot X .$$
for any $M, N\in \No$, $X\in \ca$.
\end{itemize}
\end{lema}
\pf The functor $\uhom_{\Mo}( M, -):\Mo\to \ca$ is the right adjoint of the functor $R_M:\ca\to \Mo$, $R_M(X)=X \triangleright M$. Since $R_M$ is a $\ca$-module functor then, it follows from  Lemma \ref{adjoint-module-functor} that, $\uhom_{\Mo}( M, -)$ is also a $\ca$-module functor. This implies in particular that there are natural isomorphisms 
$$\uhom_{\Mo}( M, X\triangleright N)\simeq X\ot \uhom_{\Mo}( M, N). $$
The other three isomorphisms follow in a similar way.
\epf

Let $\Mo$ be a left $\ca$-module category. There is a relation between the internal hom of $\Mo$ and $\Mo^{\op}$, stated in the next Lemma.
\begin{lema}\label{hom-op}  For any $M\in \Mo$, the functors
\begin{equation*}
 {}^{**}\uhom_\Mo(M, -), \,\uhom_{\Mo^{\op}}( -,\overline{M}):\Mo^{\bop}\to \ca,
\end{equation*}
are equivalent $\ca$-module functors. Also the functors
\begin{equation*}
\uhom_{\Mo^{\op}}( M,-)^*, \, {}^*\uhom_\Mo(-,M):\Mo\to \ca
\end{equation*} 
are equivalent $\ca$-module functors. In particular, there are natural isomorphisms
$${}^{**}\uhom_\Mo(M,N)\simeq  \uhom_{\Mo^{\op}}(\overline{N}, \overline{M}),$$
for any $M, N\in \Mo$.
\end{lema}
\pf  The functors  $D:\ca\to \ca^{\bop}$, $D(X)=X^{**},$ and $L_M :\ca^{\bop}\to \Mo^{\bop},$ $L_M(X)=X\triangleright M$, are $\ca$-module functors. A straightforward computation shows that 
$$ (L_M\circ D)^{r.a.}\simeq  \uhom_{\Mo^{\op}}( -,\overline{M}),$$
$$ D^{r.a.} \simeq {}^{**}(-),\quad (L_M)^{r.a.}\simeq  \uhom_\Mo(M,-).$$ 
Since $D$ and $L_M$ are $\ca$-module functors, then, using Lemma \ref{adjoint-module-functor} (i), it follows that, functors ${}^{**}\uhom_\Mo(M, -), \,\uhom_{\Mo^{\op}}( -,\overline{M}):\Mo^{\bop}\to \ca,$ are $\ca$-module functors. Since $ (L_M\circ D)^{r.a.}\simeq  D^{r.a.} \circ (L_M)^{r.a.}$, it follows from Lemma \ref{adjoint-module-functor} (iii) that, functors ${}^{**}\uhom_\Mo(M, -), \,\uhom_{\Mo^{\op}}( -,\overline{M}):\Mo^{\bop}\to \ca,$ are equivalent as $\ca$-module functors. The proof that,  functors $$\uhom_{\Mo^{\op}}( M,-)^*, \, {}^*\uhom_\Mo(-,M):\Mo\to \ca$$ 
are equivalent is done by showing that both functors are left adjoint of $L_M : \ca \rightarrow \Mo, L_M(X)= X \triangleright M$.
\epf

\begin{prop}\label{hom-for-algebras} Let $A\in \ca$ be an algebra. The following statements hold.
 \begin{itemize}
\item[(i)] For any $M, N\in \ca_A$, $\uhom_{\ca_A}(M,N)=(M\ot_A {}^*N)^*.$

\item[(ii)] For any $M, N\in \ca_A$, $\uhom_{(\ca_A)^{\op}}(M,N)={}^*(N\ot_A {}^*M)$.

\end{itemize}
\end{prop}
\pf Both calculations of the internal hom follow from \eqref{hom-iso-2}.
\epf 
The following  result is \cite[Lemma 3]{FS}. We include the proof since we will need later an explicit description of certain isomorphism.
\begin{lema}\label{int-hom-adjoint} Let $\Mo$ be an exact module category over $\ca$, and $F:\Mo\to \Mo$ be a $\ca$-module functor with left adjoint $F^{l.a.}:\Mo\to \Mo$. Then, there are natural isomorphisms
$$ \xi_{M,N}:\uhom(M, F(N))\to  \uhom(F^{l.a.}(M), N).$$
\end{lema}
\pf Since $F$ is a module functor, then $F^{l.a.} $ is also a module functor. Let us denote by 
$$b_{X,M}: F^{l.a.}(X\triangleright M)\to X \triangleright F^{l.a.}(M)$$
its module structure.
Let $\Omega_{M,N}: \Hom_\Mo(M, F(N))\to  \Hom_\Mo(F^{l.a.}(M), N)$ be natural isomorphisms. Take $X\in \ca$. The desired natural isomorphism is the one induced by the composition of isomorphisms
\begin{align*} &\Hom_\ca(X, \uhom(M, F(N))) \simeq \Hom_\Mo(X \triangleright M,F(N) )\\
&\simeq \Hom_\Mo(F^{l.a.}(X \triangleright M), N ) \simeq \Hom_\Mo(X \triangleright F^{l.a.}(M), N ) \simeq\\
&\simeq\Hom_\ca(X,  \uhom(F^{l.a.}(M), N).
\end{align*}     
Using isomorphisms \eqref{Hom-interno}, one can describe explicitly this isomorphism as
\begin{align}\label{uhomadj} \xi_{M,N}=\psi^Z_{F^{l.a.}(M),N}\big(\Omega_{Z \triangleright M,N}(\phi^Z_{M,F(N)}(\id_Z))b^{-1}_{Z,M} \big),
\end{align}    
where $Z=\uhom(M, F(N))$.
\epf

\subsection{ The relative Serre functor} Let $\Mo$ be a left $\ca$-module category. Following \cite{GSch}, \cite{FSS} we recall the definition of the relative Serre functor of a module category. The reader is also referred to \cite{Sh4}.
\begin{defi} A  \textit{relative Serre functor} for $\Mo$ is a pair $ (\sy_\Mo, \phi)$, where  $\sy_\Mo:\Mo\to \Mo$ is a functor equipped with natural isomorphisms 
\begin{equation}\label{Serre-equation1} \phi_{M,N}:\uhom(M,N)^*\simeq \uhom(N, \sy_\Mo(M)),
\end{equation}
for any $M, N\in \Mo$.
\end{defi}
In the next Proposition we summarize some known facts about relative Serre functors that will be used later.

\begin{prop}\label{aboute-relat-serre} Let $\Mo$ be a left module category over $\ca$. The following holds.
\begin{itemize} 
\item[(i)] $\Mo$ posses a relative Serre functor if and only if $\Mo$ is exact.

\item[(ii)]  The functor $\sy_\Mo:\Mo\to \Mo^{\bop}$ is an equivalence of $\ca$-module categories.

\item[(iii)] The natural isomorphism $\phi_{M,N}:\uhom(M,N)^*\to \uhom(N, \sy_\Mo(M)),$ is an isomorphism of $\ca$-bimodule functors.

\item[(iv)]  The relative Serre functor is unique up to isomorphism of $\ca$-module functors.
\end{itemize}\qed 
\end{prop}

\subsection{ Balanced tensor functors and Deligne tensor product}\label{Deligne-tensor}

We shall briefly recall the definition of the\textit{ relative Deligne tensor product} over a tensor category. The reader is referred to  \cite{DSS}, \cite{Gr} for more details. Assume that $\Mo$ is a right $\ca$-module category and $\No$ a left $\ca$-module category. Let $\Ac$ be a category. 

A $\ca$-\textit{balanced  functor} is a pair $(\Phi, b)$, where  $\Phi: \Mo\times \No\to \Ac$ is a functor, right exact in each variable,  equipped with natural isomorphisms $b_{M,X,N}: \Phi( M\triangleleft X, N)\to \Phi(M, X\triangleright N)$ such that it satisfies the pentagon
\begin{equation}\label{c-balanced} \Phi(\id_M, m^\No_{X,Y,N}) b_{M, X\ot Y, N}= b_{M,X,Y\rhd N} b_{M\triangleleft X, Y, N} \Phi(m^\Mo_{M,X,Y},\id_N),
\end{equation}
for any $X, Y\in \ca$, $M\in \Mo$, $N\in \No$. The natural isomorphism $b$ is called \textit{the balancing} of $\Phi$.
If $(\Phi,b),  (\widetilde{\Phi}, \widetilde{b}):  \Mo\times \No\to \Ac$ are  $\ca$-balanced functors, a $\ca$-\textit{balanced natural transformation} $\alpha: \Phi\to \widetilde{\Phi}$ is a natural transformation such that
\begin{equation}\label{c-balanced-nat} \alpha_{M, X\triangleright N} b_{M,X, N}=\widetilde{b} _{M,X, N}  \alpha_{M\triangleleft X , N},
\end{equation}
for any $X\in \ca$, $M\in \Mo$, $N\in \No$.
The \textit{balanced  tensor product}  (or sometimes called \textit{relative Deligne tensor product}) is a category $\Mo\boxtimes_\ca \No$, equipped with a $\ca$-balanced  functor $\boxtimes_\ca:  \Mo\times \No \to \Mo\boxtimes_\ca \No$ such that for any category $\Ac$ the functor
$$ \Rex(\Mo\boxtimes_\ca \No , \Ac) \to \bal( \Mo\times \No , \Ac)$$ $$ F\mapsto  F\circ \boxtimes_\ca$$
is an equivalence of categories. Here $ \bal( \Mo\times \No , \Ac)$ denotes the category of  $\ca$-balanced functors and $\ca$-balanced natural transformations.

\begin{lema}\label{balanced-comp} Let $\Mo, \widetilde{\Mo}$ be  right $\ca$-module categories and $\No, \widetilde{\No}$ be left $\ca$-module categories. If $(F,c):\widetilde{\Mo}\to \Mo$, $(G,d):\widetilde{\No}\to \No $ are right exact module functors, and  $(\Phi,b)  :  \Mo\times \No\to \Ac$ is a $\ca$-balanced functor, then $\Phi\circ (F\times G): \widetilde{\Mo} \times \widetilde{\No} \to \Ac$ is a $\ca$-balanced functor with balancing  given by
\begin{equation}\label{balcing-comp} e_{M, X, N} = \Phi(\id_{F(M)}, d^{-1}_{X,N}) b_{F(M),X,G(N)} \Phi(c_{M,X}, \id_{G(N)}),
\end{equation}
for any $M\in\widetilde{\Mo} , N\in \widetilde{\No}, X\in \ca$.
\end{lema}
\pf 
We must show that $e$ satisfies \eqref{c-balanced}. In this case we have to prove 
\begin{equation}\label{mustprove}
    \Phi (\id_{F(M)}, G(m^{\widetilde{\No}}_{X, Y, N})) e_{M, X \ot Y, N} = e_{M, X, Y \triangleright N} e_{M \triangleleft X, Y, N} \Phi (F(m^{\widetilde{\Mo}}_{M,X,Y}), \id_{G(N)}),
\end{equation}
for any $X, Y\in \ca,$ $M\in \widetilde{\Mo} , N\in \widetilde{\No}$. The left hand side of \eqref{mustprove}  is equal to
\begin{align*}
    & = \Phi (\id_{F(M)}, G(m^{\widetilde{\No}}_{X, Y, N}) d^{-1}_{X \ot Y, N}) b_{F(M), X \ot Y, G(N)} \Phi(c_{M, X \ot Y}, \id_{G(N)})\\
    & = \Phi(\id_{F(M)}, d^{-1}_{X, Y \triangleright N} (\id_X \triangleright d^{-1}_{Y,N}) m^{\No}_{X,Y,G(N)})b_{F(M), X \ot Y, G(N)}\\ & \Phi(c_{M, X \ot Y}, \id_{G(N)})\\
    & = \Phi(\id_{F(M)}, d^{-1}_{X, Y\triangleright N}) \Phi(\id_{F(M)}, \id_X \triangleright d^{-1}_{Y,N}) b_{F(M),X, Y\triangleright G(N)} b_{F(M) \triangleleft X, Y, G(N)}\\ & \Phi(m^{\Mo}_{F(M), X, Y}, \id_{G(N)}) \Phi(c_{M, X \ot Y}, \id_{G(N)}).
\end{align*}
The first equality is by the definition of $e$. The second equality is a consequence of $(G,d)$ being a  module functor, and the last equality is because $(\Phi, b)$ is a $\ca$-balanced functor. The right hand side of \eqref{mustprove}  is equal to
\begin{align*}
& = e_{M, X, Y \triangleright N} \Phi(\id_{F(M \triangleleft X)}, d^{-1}_{Y,N}) b_{F(M \triangleleft X), Y, G(N)} \Phi(c_{M \triangleleft X,Y} F(m^{\widetilde{\Mo}}_{M, X, Y}), \id_{G(N)})\\
& = e_{M, X, Y \triangleright N} \Phi(\id_{F(M \triangleleft X)}, d^{-1}_{Y,N}) b_{F(M \triangleleft X), Y, G(N)} \Phi((c^{-1}_{M,X} \triangleleft \id_Y)\\ & m^{\Mo}_{F(M), X, Y} c_{M, X \ot Y}, \id_{G(N)})\\
& = \Phi(\id_{F(M)}, d^{-1}_{X, Y \triangleright N}) b_{F(M), X, G(Y \triangleright N)} \Phi(c_{M,X}, \id_{G(Y\triangleright N)})  \Phi(\id_{F(M \triangleleft X)}, d^{-1}_{Y,N})\\
& b_{F(M \triangleleft X), Y, G(N)} \Phi((c^{-1}_{M,X} \triangleleft \id_Y) m^{\Mo}_{F(M), X, Y} c_{M, X \ot Y}, \id_{G(N)})\\
& = \Phi(\id_{F(M)}, d^{-1}_{X, Y \triangleright N}) b_{F(M), X, G(Y \triangleright N)} \Phi(c_{M,X}, \id_{G(Y \triangleright N)}) \Phi(\id_{F(M \triangleleft X)}, d^{-1}_{Y,N})\\
& \Phi(c^{-1}_{M,X}, \id_{Y \triangleright G(N)}) b_{F(M) \triangleleft X, Y, G(N)} \Phi(m^{\Mo}_{F(M), X,Y} c_{M, X \ot Y}, \id_{G(N)})\\
& = \Phi(\id_{F(M)}, d^{-1}_{X, Y \triangleright N}) b_{F(M), X, G(Y \triangleright N)} \Phi(\id_{F(M \triangleleft X)}, d^{-1}_{Y,N}) b_{F(M) \triangleleft X, Y, G(N)}\\
& \Phi(m^{\Mo}_{F(M), X,Y} c_{M, X \ot Y}, \id_{G(N)})\\
& = \Phi(\id_{F(M)}, d^{-1}_{X, Y \triangleright N}) \Phi (\id_{F(M)}, \id_X \triangleright d^{-1}_{Y,N}) b_{F(M), X, Y \triangleright G(N)} b_{F(M) \triangleleft X, Y, G(N)}\\
& \Phi(m^{\Mo}_{F(M), X,Y} c_{M, X \ot Y}, \id_{G(N)}).
\end{align*}
 The first and third equalities follow by the definition of $e$. The second equality follows since $(F,c)$ is a  module functor. The fourth equality is consequence of the naturality of $b$ for $c_{M,X}$, and the sixth equality is the naturality of $b$ for $d_{Y,N}$.
Since both sides are equal, we get the result.
\eq
\epf

The next result is well-known. 
\begin{prop}\label{basic-stuff} Let $A, B\in \ca$ be algebras. Thus, the categories $\ca_A, \ca_B$ are left $\ca$-module categories. The following assertions hold.
\begin{itemize}
\item[(i)] The functor ${}^*( -): (\ca_A)^{\op} \to {}_A\ca$ is an equivalence of right $\ca$-module categories.

\item[(ii)] The restriction of the tensor product $\ot: {}_A\ca\times \ca_B\to  {}_A\ca_B$ is a $\ca$-balanced functor, and induces an  equivalence of categories $\widehat{\ot}: {}_A\ca\boxtimes_{\ca} \ca_B \to  {}_A\ca_B$, such that
$$ \widehat{\ot}\circ \boxtimes_\ca \simeq \ot,$$
as $\ca$-balanced functors.
\item[(iii)] Assume that $\ca_A$ is an exact module category. The functor $R: {}_A\ca_B\to\Fun_\ca(\ca_A,\ca_B)$, $V\mapsto  -\ot_A V$ is an equivalence of categories.
\end{itemize}

\end{prop}
\pf 
(i) The duality functor ${}^*( -): (\ca_A)^{\op} \to {}_A\ca$ has structure of module functor with isomorphisms given by
$$\phi^l_{X^*, M}: {}^*(X^*\ot M)\to {}^*M\ot X,$$
for any $X\in \ca$, $M\in \ca_A$. Here $\phi^l$ is the natural isomorphisms described in \eqref{duals-tensor-product}. Note that we are omitting the canonical natural isomorphism ${}^*(X^*)\simeq X$.
For (ii) see \cite{DSS}. The proof of  (iii) can be found for example in \cite[Prop. 3.3]{MM}. Exactness of the module category $\ca_A$ implies that $\ot_A$ is biexact. This fact was used in \cite{MM} to prove that $R$ is a category equivalence.
\epf

\section{The (co)end for module categories}\label{Section:mcoends}

Let $\ca$ be a  tensor category and $\Mo$ be a left $\ca$-module category. Assume that $\Ac$ is a category and $S:\Mo^{\op}\times \Mo\to \Ac$ a functor equipped with natural isomorphisms
\begin{equation} \beta^X_{M,N}: S(M,X\triangleright N)\to S(X^*\triangleright M,N),
\end{equation}
for any $X\in \ca, M,N\in \Mo$. We shall say that $\beta$ is a \textit{pre-balancing} of the functor $S$.

\begin{defi} The \textit{ module end} of the pair $(S,\beta)$ is an object $E\in \Ac$ equipped with dinatural transformations $\pi_M: E\xrightarrow{ . .} S(M,M)$ such that 
\begin{equation}\label{dinat:module} S(\ev_X\triangleright \id_M,  \id_M) \pi_M= S(m_{X^*,X,M},  \id_M) \beta^X_{X\triangleright M,M} \pi_{X\triangleright M},
\end{equation}
for any $X\in \ca, M\in \Mo$, universal with this property. This means that, if $\widetilde{E}\in \Ac$ is another object with dinatural transformations $\xi_M:\widetilde{E}\xrightarrow{ . .} S(M,M)$, such that they verify \eqref{dinat:module}, there exists a unique morphism $h:\widetilde{E}\to E$ such that $\xi_M= \pi_M\circ h $.
\end{defi}
Sometimes we will denote the module end  as $\oint_{M\in \Mo} (S,\beta)$, or simply as $\oint_{M\in \Mo} S$, when the pre-balancing $\beta$ is understood from the context.
The \textit{ module coend} of the pair $(S,\beta)$ is defined dually. This  is an object $C\in \Ac$ equipped with dinatural transformations $\pi_M:S(M,M)\xrightarrow{ . .} C$ such that
\begin{equation}\label{dinat:c:module} \pi_M=\pi_{X^*\triangleright M}\beta^X_{M,X^*\triangleright M} S(\id_M, m_{X,X^*,M})S(\id_M,\coev_X\triangleright\id_M),
\end{equation}
for any $X\in \ca, M\in \Mo$, universal with this property. This means that, if $\widetilde{C}\in \Ac$ is another object with dinatural transformations  $\lambda_M: S(M,M) \xrightarrow{ . .} \widetilde{C}$ such that they satisfy \eqref{dinat:c:module}, there exists a unique morphism $g:C\to \widetilde{C}$ such that $g \circ \pi_M=\lambda_M$. The module coend will be denoted  $\oint^{M\in \Mo} (S,\beta)$, or simply as $\oint^{M\in \Mo} S$.

\medbreak

A similar definition can be made for \textit{right} $\ca$-module categories. Let $\Bc$ be a category, and $\No$ be a right $\ca$-module category endowed with a functor $S:\No^{\op}\times \No\to \Bc$ with a \textit{pre-balancing}
$$ \gamma^X_{M,N}: S(M\triangleleft X, N)\to S(M, N\triangleleft {}^* X),$$
for any $M,N\in \No$, $X\in \ca$. 
\begin{defi} The \textit{module end} for $S$ is an object $E\in \Bc$ equipped with dinatural transformations  $\lambda_N:E\xrightarrow{ .. } S(N,N)$ such that 
\begin{equation}\label{dinat:module-left}
\lambda_N= S(\id_N, \id_N \triangleleft \ev_X)S(\id_N, m^{-1}_{N,X,{}^* X}) \gamma^X_{N,N\triangleleft X} \lambda_{N\triangleleft X},
\end{equation}
for any $N\in \No$, $X\in \ca$. We shall also denote this module end by $\oint_{N\in \No} (S, \gamma)$.

Similarly, the \textit{module coend} is an object $C\in \Bc$ with dinatural transformations $\lambda_N: S(N,N)\xrightarrow{ .. } C$  such that
\begin{equation}\label{dinat:c-module-left}
\lambda_N S(\id_N  \triangleleft \coev_X,\id_N)= \lambda_{N  \triangleleft  {}^* X} \gamma^X_{N  \triangleleft  {}^* X, N}S(m^{-1}_{N,   {}^* X, X}, \id_N),
\end{equation}
for any $N\in \No$, $X\in \ca$. We shall also denote this module coend by $\oint^{N\in \No} (S, \gamma)$.
\end{defi}

In the next Proposition we collect some  results about module ends that  generalize   well-known results in the theory of (co)ends. The proofs follow the same lines as the ones in usual ends. For the sake of completeness we include them.

\begin{prop}\label{equivalences-mod-end} Assume that $\Mo, \No$ are left $\ca$-module categories, and $S, \widetilde{S}:\Mo^{\op}\times \Mo\to \Ac$ are functors equipped with pre-balancings $$\beta^X_{M,N}: S(M,X\triangleright N)\to S(X^*\triangleright M,N),$$ $$ \widetilde{\beta}^X_{M,N}: \widetilde{S}(M,X\triangleright N)\to \widetilde{S}(X^*\triangleright M,N),$$ $X\in \ca, M,N\in \Mo$. The following assertions holds
\begin{itemize}
\item[(i)] Assume that the module ends $\oint_{M\in \Mo} (S,\beta), \oint_{M\in \Mo} (\widetilde{S},\widetilde{\beta})$ exist and have dinatural transformations $\pi, \widetilde{\pi}$, respectively. If  $\gamma:S\to \widetilde{S}$ is a natural transformation  such that 
\begin{equation}\label{gamma1}
 \widetilde{\beta}^X_{M,N} \gamma_{(M,X\triangleright N)}=\gamma_{(X^*\triangleright M,N)} \beta^X_{M,N},
\end{equation}
then there exists a unique map $\widehat{\gamma}: \oint_{M\in \Mo} (S,\beta)\to \oint_{M\in \Mo} (\widetilde{S},\widetilde{\beta})$ such that   $\widetilde{\pi}_M \widehat{\gamma}= \gamma_{(M,M)} \pi_M$ for any $M\in \Mo$. If $\gamma$ is a natural isomorphism, then  $\widehat{\gamma}$ is an isomorphism.
\item[(ii)] If the end  $\oint_{M\in \Mo} (S,\beta)$ exists, then for any object $U\in \Ac$,  the end $\oint_{M\in \Mo} \Hom_\Ac(U, S( -, -))$ exists, and there is an isomorphism
$$\oint_{M\in \Mo} \Hom_\Ac(U, S( -, -))\simeq  \Hom_\Ac(U, \oint_{M\in \Mo} (S,\beta) ).$$
Moreover, if $\oint_{M\in \Mo} \Hom_\Ac(U, S( -, -))$ exists for any $U\in \Ac$, then the end $\oint_{M\in \Mo} (S,\beta)$ exists.
\item[(iii)] Assume $F:\Ac\to \Ac'$ is a right exact functor. Then, there is an isomorphism
$$ F(\oint_{N\in \No} (S,\beta) )\simeq \oint_{N\in \No} (F\circ S,F(\beta)).$$

\item[(iv)]  If $F: \Mo\to \No$ is an equivalence of $\ca$-module categories, then there is an isomorphism
$$\oint_{N\in \No} S \simeq \oint_{M\in \Mo} S(F(-),F(-)).$$
\end{itemize}

\end{prop}

\pf   (i). For any $M\in \No$ define 
$\lambda_M: \oint_{N\in \No} (S,\beta)\to \widetilde{S}(M,M)$ as $\lambda_M= \gamma_{(M,M)} \pi_M$. It follow straightforward that $\lambda$ is dinatural and  since $\gamma$ satisfies \eqref{gamma1}, then $\lambda$ satisfies \eqref{dinat:module}. By the universality of the module end, there exists a morphism $\widehat{\gamma}: \oint_{M\in \Mo} (S,\beta)\to \oint_{M\in \Mo} (\widetilde{S},\widetilde{\beta})$ such that   $\widetilde{\pi}_M \widehat{\gamma}= \lambda_M= \gamma_{(M,M)} \pi_M$.

(ii). Let us assume that $\oint_{M\in \Mo} (S,\beta)$ exists, and has associated to it dinatural transformations $\pi_N: \oint_{M\in \Mo} (S,\beta)\to S(N,N)$. For any $U\in \Ac$, the pre-balancing for the functor $\Hom_\Ac(U, S( -,-)) $ is defined as
$$\beta^U_{X,M,N}: \Hom_\Ac(U, S( M, X\triangleright N))\to \Hom_\Ac(U, S( X^*\triangleright M, N)),$$
$$\beta^U_{X,M,N}(f)=\beta^X_{M,N}\circ f. $$
Also define 
$$\pi^U_N:  \Hom_\Ac(U, \oint_{M\in \Mo} (S,\beta) ) \to   \Hom_\Ac(U, S(N,N)),$$
$$ \pi^U_N(f)=\pi_N\circ f.$$
It follows by a straightforward computation that, $\pi^U$ is a dinatural transformation, and they satisfy \eqref{dinat:module} using $\beta^U$. It also follows easily that $\Hom_\Ac(U, \oint_{M\in \Mo} (S,\beta) ) $ together with $\pi^U$ satisfy the universal property of the module end, thus 
$$\oint_{M\in \Mo} \Hom_\Ac(U, S( -, -))\simeq \Hom_\Ac(U, \oint_{M\in \Mo} (S,\beta) ).  $$
Now, let us assume that $\oint_{M\in \Mo} \Hom_\Ac(U, S( -, -))$ exists for any $U\in \Ac$.  Using item (i),   we can define a functor 
$$ F:\Ac^{\op}\to \vect_\ku,$$ $$F(U)=\oint_{M\in \Mo} \Hom_\Ac(U, S( -, -)).$$
We shall prove that $F$ es left exact, and thus it is representable. The object representing the functor $F$ will be a candidate for the module end $\oint_{M\in \Mo} (S,\beta)$.

 For any $M\in \Mo$, and any $f:U\to V$ in $\Ac$, denote  $$(\alpha_f)_M:\Hom_\Ac(V, S( M, M)) \to  \Hom_\Ac(U, S( M, M))$$  $$(\alpha_f)_M(g)=g\circ f.$$ To prove that $F$ is left exact, we need to show that, for any morphism $f:U\to V$ in $\Ac$, $F(\cok(f))=\Ker (F(f))$. Let be $q=\cok(f):V\to C$, and $l:K\to F(V)$ be a $\ku$-linear map such that $F(f)\circ l=0$. Then 
 $$(\alpha_f)_M \circ \pi^V_M\circ l=\pi^U_M\circ F(f)\circ l=0. $$
The second equality follows from item (i). Since $\ker(\alpha_f)=\alpha_q$, there exists a map
$$h_M:K\to \Hom_\Ac(C,S( M, M)))$$
such that $(\alpha_q)_M\circ h_M=\pi^V_M\circ l$. It is not difficult to prove that $h$ is a dinatural transformation, and they satisfy \eqref{dinat:module}  (using the isomorphisms $\beta^C$). By the universal property of the module end, there exists a morphism $\phi:K\to F(C)$ such that $h_M=\pi^C_M\circ \phi$. It follows from item (i) that
$$ (\alpha_q)_M\circ \pi^C_M\circ \phi=\pi^V_M\circ F(q)\circ \phi.$$
But also
$$ (\alpha_q)_M\circ \pi^C_M\circ \phi= (\alpha_q)_M \circ h_M=\pi^V_M\circ l,$$
whence $l= F(q)\circ \phi$ and therefore $ F(q)=\ker(F(f))$. Hence $F$ is represented by an object $E\in \Ac$; $F(U)=\Hom_\Ac(U, E)$. The maps $\delta_M:E\to S(M,M),$ $\delta_M= \pi^E_M(\id_E)$ are dinatural transformations, and they satisfy \eqref{dinat:module}. It follows by a straightforward computation that $E$ together with $\delta$ satisfy the universal property of the module end, thus $E\simeq \oint_{M\in \Mo} (S,\beta)$.
\medbreak

The proof of (iii) and (iv) is straightforward.
\epf

\begin{rmk} Off course that, similar results to those presented in Proposition \ref{equivalences-mod-end} can be stated for module coends, and also for module (co)ends for right module categories.
\end{rmk}

\subsection{Relation between module (co)ends for right and left module categories}

Let $\Ac$ be a category. Let $\Mo$ be a left $\ca$-module category, and a functor $S:\Mo^{\op}\times \Mo\to \Ac$ equipped with a pre-balancing $\beta^X_{M,N}: S(M,X\triangleright N)\to S(X^*\triangleright M,N)$. Then $\No=\Mo^{\op}$ is a right $\ca$-module category. We can consider the functor 
$$ S^{\op}:\No^{\op}\times \No\to \Ac^{\op}.$$
It posses a pre-balancing 
$$ \gamma^X_{M,N}: S^{\op}(M\triangleleft X, N)\to S^{\op}(M, N\triangleleft {}^* X),$$
$$\gamma^X_{M,N}= \beta^X_{M,N}. $$
Note that, the pre-balancing $\gamma$ is considered as a morphism in $\Ac^{\op}$. The next result is straightforward.

\begin{lema}\label{left-right} There are isomorphisms
$$\overline{\oint_{M\in \Mo} (S,\beta)}\simeq \oint^{\overline{M}\in \No} (S^{\op},\gamma), $$
$$\overline{\oint^{M\in \Mo} (S,\beta)}\simeq \oint_{\overline{M}\in \No} (S^{\op},\gamma). $$\qed
\end{lema}

A similar result also holds starting from a right $\ca$-module category $\No$, and a functor $T:\No^{\op}\times \No\to \Ac$ equipped with a pre-balancing $$\gamma^X_{M,N}: T(M\triangleleft X, N)\to T(M, N\triangleleft {}^* X).$$
If $\Mo=\No^{\op}$, then $\Mo$ is a left $\ca$-module category, and we can consider the functor 
$$ T^{\op}: \Mo^{\op}\times \Mo\to \Ac^{\op}$$
together with a pre-balancing
$$ \beta^X_{M,N}: T^{\op}(M,X\triangleright N)\to T^{\op}(X^*\triangleright M,N)$$
$$\beta^X_{M,N}=\gamma^{X^{**}}_{M,N}. $$
The next result is a straightforward consequence of the definitions of module (co)end.
\begin{lema}\label{right-left} There are isomorphisms
$$\overline{\oint_{N\in \No} (T,\gamma)}\simeq \oint^{\overline{M}\in \Mo} (T^{\op},\beta), $$
$$\overline{\oint^{M\in \No} (T,\gamma)}\simeq \oint_{\overline{M}\in \Mo} (T^{\op},\beta). $$\qed
\end{lema}

\subsection{Parameter theorem for module ends}\label{parameter-coends}

Let $\ca$ be a  tensor category and $\Mo$ a left $\ca$-module category. Also, let $\Ac, \Bc$ be categories. We start with a functor $S:\Mo^{\op}\times \Mo\to \Fun(\Ac, \Bc)$ equipped with pre-balancing 
 $\beta^X_{M,N}: S(M,X\triangleright N)\to S(X^*\triangleright M,N),$
for any $X\in \ca, M,N\in \Mo$. If the end  $\oint_{M\in \Mo} (S,\beta)$ exists, it is an object in the category $\Fun(\Ac, \Bc)$; we denote this functor as
$$ \big(\oint_{M\in \Mo} (S,\beta) \big)(-):\Ac\to \Bc.$$

Alternatively, we can do the following construction. For any $A\in \Ac$ we get a functor $S_A: \Mo^{\op}\times \Mo\to \Bc$, $S_A(M,N)=S(M,N)(A)$. This functor comes with a pre-balancing 
$$ \beta^A_{X,M,N}: S_A(M,X\triangleright N)\to S_A(X^*\triangleright M,N),$$
$$ \beta^A_{X,M,N}= (\beta^X_{M,N})_A,$$
for any $X\in \ca, M,N\in \Mo$. If the module end $ \oint_{M\in \Mo} (S_A,\beta^A)$ exists, it is an object in $\Bc$, and it defines a functor $\mathfrak{S}: \Ac\to \Bc.$ The proof of the next result follow straightforward.

\begin{teo}\label{parameter-end} Provided all ends $ \oint_{M\in \Mo} (S_A,\beta^A)$ exist, the functor $\mathfrak{S}$ has a canonical structure of module end for the functor $S$. We write
$$ \mathfrak{S}= \big(\oint_{M\in \Mo} (S,\beta) \big)(-).$$\qed
\end{teo}

\begin{rmk} Similar results can be obtained for module coends, and also for right $\ca$-module categories.
\end{rmk}

\subsection{Restriction of module (co)ends to tensor subcategories}
In this Section, we shall show that the module (co)end coincides with the usual (co)end in the case the tensor category is $\vect_\ku$. We also study what happens with the module (co)end when we restrict to a tensor subcategory.
\smallbreak

Let $\ca$ be a tensor category and $\Do \subseteq \ca$ be a tensor subcategory of $\ca$. Assume also that $\Mo$ is a left $\ca$-module category.  We can consider the restricted $\Do$-module category $\Res^\Do_\ca \Mo$.  The next result is a straightforward consequence of the definition of module (co)ends.
 \begin{prop}\label{restriction-sucategory}
  Let $S:\Mo^{\op}\times \Mo\to \Ac$ be a functor equipped with pre-balancing $\beta^X_{M,N}: S(M,X\triangleright N)\to S(X^*\triangleright M,N)$. 
 \begin{itemize}
\item[(i)]  There exists a monomorphism in $\Ac$
$$ \oint_{M\in \Mo} (S,\beta) \hookrightarrow\oint_{M\in \Res^\ca_\Do \Mo} (S,\beta).$$

\item[(ii)]  There exists an epimorphism in $\Ac$ $$ \oint^{M\in \Res^\Do_\ca \Mo} (S,\beta) \twoheadrightarrow \oint^{M\in \Mo} (S,\beta) .$$
 \end{itemize}
  
 \qed
 \end{prop}

\begin{rmk} Similar result  obtained in Proposition \ref{restriction-sucategory},  is valid for right module categories. 
\end{rmk}

The next result says that, the module (co)end coincides with the usual one in the case $\ca=\vect_\ku$.
 \begin{prop}\label{restriction-vect} Let $\Mo, \Ac$ be categories, and $S: \Mo^{\op}\times \Mo\to \Ac$ be a functor. In particular $\Mo$ is a left $\vect_\ku$-module category. The functor $S$ has a canonical pre-balancing $\beta$ such that there are isomorphisms
 $$\int_{M\in \Mo} S\simeq \oint_{M\in \Mo} (S,\beta),$$
 $$\int^{M\in \Mo} S\simeq \oint^{M\in \Mo} (S,\beta).$$
 \end{prop}
 \pf We shall prove the first isomorphism concerning the usual end and the module end. The other isomorphism for the coend follows similarly. For this, we will show that, for such a functor $S$ there exists a canonical pre-balancing $\beta$ such that \textit{any} dinatural transformation $\pi_M:E\xrightarrow{ .. }S(M,M)$ satisfies \eqref{dinat:module}.

 Since $\Mo$ is a finite abelian $\ku$-linear category, there exists a finite dimensional $\ku$-algebra $A$ such that, $\Mo$ is equivalent to the category of finite dimensional right $A$-modules $\moda_A$. The action of $\vect_\ku$ on $\moda_A$ is
 $$\triangleright:\vect_\ku\times \moda_A \to \moda_A, $$
 $$X\triangleright M=X\otk M, $$
 for any $X\in \vect_\ku$, $M\in  \moda_A$. The right action of $A$ on  $X\otk M$ is on the second tensorand. For  any $X, Y\in \vect_\ku$, $M\in  \moda_A$  the associativity of this module category is
 $$m_{X,Y,M}:(X\otk Y)\otk M\to X\otk (Y\otk M),$$
 $$m_{X,Y,M}((x\ot y)\ot m)=x\ot (y\ot m). $$
 
 For any $X\in \vect_\ku$, $x\in X$, we  denote by $\delta_x:X\to \ku$ the unique linear transformation that sends $x$ to 1, and any element of a direct complement of $<x>$ to 0. If $M\in  \moda_A$, $X\in \vect_\ku$, $x\in X$ we shall denote by
 $$\iota^M_x:M\to X\triangleright M,\quad p^M_x:X\triangleright M\to M, $$
 $$\iota^M_x(m)=x\ot m, \quad  p^M_x(y\ot m)=\delta_x(y) m, $$
 for any $y\in X, m\in M$. Let $(x_i), (f_i)$ be a pair of dual basis of $X$ and $X^*$ respectively. For any $x\in X, f\in X^*$ it is easy to verify that
 $$ \sum_i \delta_{x_i} (x)  \delta_{f_i} (f) =f(x).$$
 This equality implies that
 \begin{equation}\label{dual-basis-ev} \ev_X\ot \id_M= \sum_i  p^M_{x_i} p^{X\otk M}_{f_i} m_{X^*,X,M}.
 \end{equation}
 Also, one can verify that 
 \begin{equation}\label{incl-proj} \sum_i\, \iota^M_{x_i} p^M_{x_i} =\id_{X\otk M}, \quad 
 p^M_{x} \iota^M_{y}=\delta_x(y) \id_M.
 \end{equation}
 
 For any $M, N\in  \moda_A$ let us denote
$$\beta^X_{M,N}: S(M,X\triangleright N)\to S(X^*\triangleright M,N),$$ 
$$\beta^X_{M,N}=\oplus_{i} S(p^M_{f_i}, p^N_{x_i}),$$ 
where $(x_i), (f_i)$ is a pair of dual basis of $X$ and $X^*$, respectively.  Using \eqref{incl-proj}, one can check that, $\beta^X_{M,N}$ is an isomorphism with inverse 
$$ \oplus_{i} S(\iota^M_{f_i}, \iota^N_{x_i}): S(X^*\triangleright M,N)\to S(M,X\triangleright N).$$
 Let $E\in \Ac$ be an object and $\pi_M:E\xrightarrow{ . .} S(M,M)$  a dinatural transformation. Let us show that, $\pi$ satisfies equation \eqref{dinat:module}. Let $X\in \vect_\ku$, $M\in  \moda_A$ and let   $(x_i), (f_i)$ be a pair of dual basis of $X$ and $X^*$. The right hand side of equation \eqref{dinat:module} is
\begin{align*} S(m_{X^*,X,M},  \id_M)& \beta^X_{X\triangleright M,M} \pi_{X\triangleright M}=\oplus_{i}  S(m_{X^*,X,M},  \id_M) S(p^{X\triangleright M}_{f_i}, p^M_{x_i}) \pi_{X\triangleright M}\\
&=\oplus_{i}  S(m_{X^*,X,M},  \id_M)  S( p^M_{x_i}p^{X\triangleright M}_{f_i},\id_M) \pi_M\\
&=\oplus_{i}  S( p^M_{x_i}p^{X\triangleright M}_{f_i}m_{X^*,X,M},\id_M) \pi_M\\
&= S(\ev_X\ot\id_M, \id_M) \pi_M.
\end{align*}
The second equality follows from the dinaturality of $\pi$, and the last equality follows from \eqref{dual-basis-ev}.
 \epf

Combining Proposition \ref{restriction-sucategory} and Proposition \ref{restriction-vect} we have the next result.

\begin{cor} Assume $\Mo$ is a left $\ca$-module category, $\Ac$ is a category, and  $S: \Mo^{\op}\times \Mo\to \Ac$ is a functor. Let  $\lambda_M: \int_{M\in  \Mo} S\xrightarrow{ . .}  S(M,M)$ be the associated dinatural transformation of the usual end. There exists a monomorphism
$$\varphi: \oint_{M\in \Mo} (S,\beta) \to \int_{M\in  \Mo} S $$
such that 
\begin{equation}\label{comparison-1} S(\ev_X\triangleright \id_M,  \id_M) \lambda_M \varphi= S(m_{X^*,X,M},  \id_M) \beta^X_{X\triangleright M,M} \lambda_{X\triangleright M} \varphi,
\end{equation} 
for any $X\in \ca, M\in \Mo$.\qed
\end{cor}

Using the above Corollary, we can give another characterization of the module end. Essentially, this new description says that the module end is a subobject of the usual end, and it is universal among those subobjects with morphisms that satisfy \eqref{comparison-1}.

\begin{prop}\label{ach} Let $(E, \psi)$ be a pair, where
\begin{itemize}
\item $E$ is an object in $\Ac$;

\item $\psi: E\to \int_{M\in  \Mo} S$ is a morphism such that
\begin{equation}\label{comparison-2} S(\ev_X\triangleright \id_M,  \id_M) \lambda_M \psi= S(m_{X^*,X,M},  \id_M) \beta^X_{X\triangleright M,M} \lambda_{X\triangleright M} \psi.
\end{equation} 
\end{itemize}
Then, there exists a unique map $h: E\to \oint_{M\in \Mo} (S,\beta)$ such that $\psi=\varphi\circ h. $\qed
 \end{prop}

\section{Applications to the theory of representations of tensor categories}\label{Section:application}
Throughout this section $\ca$ will denote a  tensor category.

\subsection{Natural module transformations as an end}

For a pair of functors $F, G:\Ac\to \Bc$ between two abelian categories $\Ac, \Bc$, it is well known that, there is an isomorphism 
$$ \Nat(F,G)\simeq \int_{A\in \Ac} \Hom_\Bc(F(A), G(A)).$$
In this Section, we  generalize this result when $F$ and $G$ are $\ca$-module functors.

Let $\Mo, \No$ be $\ca$-module categories, and $(F,c), (G,d):\Mo \to \No$ be module functors. The functor $$\Hom_{\No}(F(-), G(-)): \Mo^{\op}\times \Mo\to \vect_\ku$$ evaluated on functions  $f:M\to M'$, $g:N\to N'$ in $\Mo$, is $$\Hom_{\No}(F(f), G(g))(\alpha)=G(g)\circ \alpha\circ F(f),$$ for any $\alpha\in \Hom_{\No}(F(M), G(N))$.

\begin{prop}\label{end-natural-transf} For any pair of $\ca$-module functors $(F, c), (G, d)$ the functor $$\Hom_{\No}(F(-), G(-)): \Mo^{\op}\times \Mo\to \vect_\ku$$
has a pre-balancing  given by
\begin{equation}\label{beta-for-homs}
 \beta^X_{M,N}: \Hom_{\No}(F(M), G(X\triangleright N))\to \Hom_{\No}(F(X^*\triangleright M), G(N))
\end{equation}
$$ \beta^X_{M,N}(\alpha)= (ev_X\triangleright \id_{G(N)}) m^{-1}_{X^*,X,G(N)} (\id_{X^*}\triangleright d_{X,N}\alpha)c_{X^*,M},$$
for any $X\in \ca, M, N\in \Mo$.
There is an isomorphism $$ \Nat_{\!m}(F,G)\simeq \oint_{M\in \Mo} (\Hom_{\No}(F(-), G(-)), \beta).$$
\end{prop}
\pf    It follows straightforward that, $ \beta^X_{M,N}$ are  natural isomorphisms with inverses given by
\begin{equation}\label{beta-inverse}
\big(\beta^X_{M,N}\big)^{–1}(\alpha)=d^{-1}_{X,N}(\id_X\triangleright \alpha c^{-1}_{X^*,M})m_{X,X^*,F(M)}(\coev_X\triangleright \id_{F(M)}).
\end{equation} 

For any $M\in \Mo$, define $\pi_M: \Nat_{\!m}(F,G)\to \Hom_{\No}(F(M), G(M))$ by $\pi_M(\alpha)=\alpha_M$. It follows easily that, $\pi$ is a dinatural transformation. Let us show that $\pi$ satisfies  \eqref{dinat:module}. Let $\alpha\in  \Nat_{\!m}(F,G)$, $M\in \Mo$, then the left hand side of \eqref{dinat:module} evaluated in $\alpha$ is equal to 
$$\alpha_M  F(\ev_X\triangleright \id_M).$$
The right hand side of \eqref{dinat:module} evaluated in $\alpha$ is equal to 
\begin{align*} &=(ev_X\triangleright \id_{G(M)})m^{-1}_{X^*,X,G(M)}(\id_{X^*}\triangleright d_{X,M}\alpha_{X\triangleright M}) c_{X^*,X\triangleright M} F(m_{X^*,X,M})\\
&=(ev_X\triangleright \id_{G(M)})m^{-1}_{X^*,X,G(M)} (\id_{X^*}\triangleright( \id_X\triangleright \alpha_M) c_{X,M})c_{X^*,X\triangleright M} F(m_{X^*,X,M})\\
&=(ev_X\triangleright \id_{G(M)}) (\id_{X^*\ot X}\triangleright \alpha_M) c_{X^*\ot X, M}\\
&= \alpha_M(ev_X\triangleright \id_{F(M)})c_{X^*\ot X, M}= \alpha_M  F(\ev_X\triangleright \id_M).
\end{align*}
The second equality follows since $\alpha$ is a module natural transformation and satisfies \eqref{modfunctor3},  the third equality follows by the naturality of $m$ and since $c$ satisfies \eqref{modfunctor1} and the last one follows from the naturality of $c$.

Let $E$ be a vector space equipped with a dinatural transformation $\xi_M:E\to \Hom_{\No}(F(M), G(M))$ such that \eqref{dinat:module} is satisfied. Define $h:E\to  \Nat_{\!m}(F,G)$ as follows. For any $v\in E$, $M\in \Mo$, 
$ h(v)_M=\xi_M(v)$. It is clear, by definition, that $\pi\circ h=\xi$. We must prove that, for any $v\in E$, $h(v)$ is a natural module transformation, that is, we must show that equation \eqref{modfunctor3} is fulfilled, which  in this case is
\begin{equation}\label{dinat-h1} d_{X,M} \xi_{X\triangleright M}(v)=(\id_X\triangleright \xi_M(v)) c_{X,M},
\end{equation}
for any $X\in \ca$, $M\in \Mo$. Since $\xi$ satisfies \eqref{dinat:module}, then
$$ \big(\beta^X_{X\triangleright M,M}\big)^{-1}\big(\xi_M(v) F(ev_X\triangleright \id_M) F(m^{-1}_{X^*,X,M})\big)= \xi_{X\triangleright M}(v),$$
for any $v\in E$. Using the definition of $\big(\beta^X_{X\triangleright M,M}\big)^{-1}$ given in \eqref{beta-inverse}, this equation is equivalent to
\begin{align*} & d_{X,M} \xi_{X\triangleright M}(v)=\big(\id_X\triangleright \xi_M(v)F((ev_X\triangleright \id_M )m^{-1}_{X^*,X,M})c^{-1}_{X^*, X\triangleright M} \big)\\& m_{X,X^*,F(X\triangleright M)} (\coev_X\triangleright \id_{F(X\triangleright M)})\\
&= \big(\id_X\triangleright \xi_M(v)F(ev_X\triangleright \id_M) \big) 
(\id_X\triangleright c^{-1}_{X^*\ot X, M} m^{-1}_{X^*, X, F(M)}\\ &(\id_{X^*} \triangleright c_{X,M})) m_{X,X^*,F(X\triangleright M)}  (\coev_X\triangleright \id)\\
&= \big(\id_X\triangleright \xi_M(v)\big) (\id_X\ot ev_X\triangleright\id_M)(\id_{X\ot X^*}\triangleright c_{X,M})(\coev_X\triangleright \id)\\
&=\big(\id_X\triangleright \xi_M(v)\big) c_{X,M}.
\end{align*}
The second equality follows from \eqref{modfunctor1}, the third equality follows from the naturality of $c$, and the last one follows from the rigidity of $\ca$. Hence, $\Nat_{\!m}(F,G)$ satisfies the required universal property.
\epf

\subsection{On the category of module functors}\label{Subsection:modulefunctors}

Assume that $\ca, \Ec, \Do,$ are  tensor categories. Assume also that $\Mo$ is a $(\ca,\Ec)$-bimodule category, and  that $ \No$ is a $(\ca,\Do)$-bimodule category. Then, we can consider the functors
\begin{equation}\label{L-equivalence1} L=L_{\Mo,\No}: \Mo^{\op}\boxtimes_\ca \No \to \Fun_\ca(\Mo^{\bop}, \No),\end{equation}
$$L(M\btc N)= \uhom_{\Mo^{\op}}( -, M)\triangleright N,$$
\begin{equation}\label{L-equivalence2} \widetilde{L}=\widetilde{L}_{\Mo,\No}: \Mo^{\op}\boxtimes_\ca \No \to \Fun_\ca(\Mo, \No)
\end{equation}
$$\widetilde{L}(M\btc N)= \uhom_{\Mo^{\op}}( M, -)^*\triangleright N,$$
Both functors are equivalences of $(\Ec,\Do)$-bimodule categories. 
This fact was proven in \cite[Thm. 3.20] {Gr}, see also \cite{DSS}.
The bimodule structure on the functor category $\Fun_\ca(\Mo, \No)$ is described in \eqref{actions-funct-categ}. We will give another proof of the fact that $L$ and $\widetilde{L} $ are category equivalences, and we shall  show an explicit description of a quasi-inverse  using the module end of some functor in an analogous way as \cite[Lemma 3.5]{Sh2}.

\begin{rmk} Our choice of  the definition of the action on the right module category $\Mo^{\op}$, given in Section \ref{bimodule-categories}, using the right dual, is justified by the definition of the functors $L$ and $\widetilde{L}$. If one changes the action on $\Mo^{\op}$ using left duals, then one has to modify the definition of the functors $L$ and $\widetilde{L}$, so that they are well-defined.
\end{rmk}

For later use, let us explain explicitly what it means that $ \widetilde{L}$ is a bimodule functor. For any $Z\in \Do$, $W\in \Ec$, $M\in \Mo, N\in \No$ we have natural isomorphisms
\begin{equation}\label{L-is-bimodule1}
 \widetilde{L}(W\triangleright M\boxtimes_\ca N)\simeq   \widetilde{L}(M\boxtimes_\ca N)\circ ( - \triangleleft W), 
\end{equation}
\begin{equation}\label{L-is-bimodule2}
 \widetilde{L}( M\boxtimes_\ca N\triangleleft Z)\simeq (-\triangleleft Z)\circ \widetilde{L}(M\boxtimes_\ca N).
\end{equation}

Assume that,  $\Mo, \No$ are exact indecomposable as left $\ca$-module categories, then there exist algebras $A, B\in \ca$ such that $\Mo\simeq \ca_A$, $\No\simeq\ca_B$ as module categories.  Recall that, if $M\in\ca_A$ then, by Lemma \ref{dual-module1} (i),   ${}^*M$ has structure of left $A$-module. Exactness of the module category $\Mo$ is needed to use Proposition \ref{basic-stuff} (iii).

\begin{lema}\label{aboutL} Assume  as above that $\Mo=\ca_A$, $\No=\ca_B$. Denote by $(\sy_\Mo, \phi)$ a relative Serre functor associated to $\Mo$. Then, the following statements hold.
\begin{itemize}
\item[(i)] The functor $\widetilde{L}_{\Mo,\No}: \Mo^{\op}\boxtimes_\ca \No \to \Fun_\ca(\Mo, \No)$ is equivalent to the composition 
 of functors
$$ (\ca_A)^{\op}\boxtimes_{\ca} \ca_B \xrightarrow{^*( -) \boxtimes\Id} {}_A\ca\boxtimes_{\ca} \ca_B \xrightarrow{\ot}  {}_A\ca_B \xrightarrow{ R } \Fun_\ca(\ca_A,\ca_B).$$
   Recall the definition of the functor $R$ given in Proposition \ref{basic-stuff}, $R:   {}_A\ca_B\to \Fun_\ca(\ca_A,\ca_B),$  $R(V)(X)=X\ot_A V$. In particular, it follows that $\widetilde{L}$ is a category equivalence.
  \item[(ii)] For any $M\in \Mo$, $N\in \No$, there exists a natural isomorphism of module functors 
 \begin{equation}\label{eles-serre}
 \widetilde{L}_{\Mo,\No}(M\btc N)\simeq L_{\Mo,\No}(M\btc N)\circ \sy_{\Mo}.
 \end{equation} 
 In particular $L$ is also an equivalence of categories.
  \end{itemize}
\end{lema}
\pf Existence of the relative Serre functor $\sy_\Mo$ is ensured by the fact that $\Mo$ is exact. Part (i) follows from the computation of the internal hom given in Proposition \ref{hom-for-algebras} (ii). Let us prove (ii).
It follows from Lemma \ref{hom-op} that, functors 
$$\uhom_{\Mo^{\op}}(M, -)^*, \, {}^*\uhom_{\Mo}(-, M): \Mo\to \ca, $$
are equivalent as $\ca$-module functors. Also, it follows from  Lemma \ref{hom-op} that,  functors 
$${}^{**}\uhom_\Mo(M, \sy_{\Mo}(-)), \,\uhom_{\Mo^{\op}}( \sy_{\Mo}(-),\overline{M}):\Mo\to \ca $$
are equivalent as $\ca$-module functors. This implies that, $\widetilde{L}_{\Mo,\No}(M\btc N)$ is equivalent to the $\ca$-module functor
$${}^*\uhom_{\Mo}(-, M)\triangleright N: \Mo\to \No, $$
and $L_{\Mo,\No}(M\btc N)\circ \sy_{\Mo}$ is equivalent to the $\ca$-module functor
$${}^{**}\uhom_\Mo(M, \sy_{\Mo}(-)) \triangleright N: \Mo\to \No. $$
The natural isomorphisms $\phi_{U,V}:\uhom(U,V)^*\to \uhom(V, \sy_\Mo(U))$ induce an isomorphism of $\ca$-module functors
$${}^{**}\phi_{-,M}\triangleright \id_N: {}^*\uhom_{\Mo}(-, M)\triangleright N \to {}^{**}\uhom_\Mo(M, \sy_{\Mo}(-))\triangleright N.$$
And this finishes the proof of the Lemma.
\epf

In what follows, we shall give an explicit description of a quasi-inverse of the functor $\widetilde{L}$ using the module end. 
For any module functor  $(F, c)\in \Fun_\ca(\ca_A,\ca_B)$ we  introduce some auxiliary functors $S_F, D_F, \ele_F, \ere_F$ that, later, we will compute its module end. 
\medbreak

Define 
\begin{equation}\label{defin-sf} \begin{split}
S_F:(\ca_A)^{\op}\times \ca_A\to (\ca_A)^{\op}\boxtimes_\ca \ca_B,\\
S_F(M,N)=\overline{M}\boxtimes_\ca F(N),
\end{split}
\end{equation}
endowed with 
a pre-balancing
 $$\beta^X_{M,N}: S_F(M,X\triangleright N)\to S_F(X^*\triangleright M,N)$$
 $$\beta^X_{M,N}= b^{-1}_{M,X,N} (\id_{M} \boxtimes_\ca  c_{X,N}).$$ 
 Also 
 \begin{equation}\label{defin-df} \begin{split} D_F :(\ca_A)^{\op}\times \ca_A\to (\ca_B)^{\op}\boxtimes_\ca \ca_A,\\
 D_F(M,N)=\overline{F(M)}\boxtimes_\ca N,
\end{split}
\end{equation}
endowed  with 
a pre-balancing
 $$\delta^X_{M,N}: D_F(M,X\triangleright N)\to D_F(X^*\triangleright M,N)$$
 $$\delta^X_{M,N}= (c^{-1}_{X^*, M}\btc \id_N) b^{-1}_{M,X,N}.$$
 Here  $b_{M,X,N}: X^*\triangleright M\boxtimes_\ca N \to M\boxtimes_\ca X\triangleright N$ is the balancing associated to the Deligne tensor product $\boxtimes_\ca$, see Section \ref{Deligne-tensor}.

 We also have  functors 
 \begin{equation}\label{defin-lrf} \begin{split} \ele_F, \ere_F:(\ca_A)^{\op}\times \ca_A\to {}_A\ca_B,\\ 
 \ere_F(\overline{M},N)={}^*M\ot F(N), \\ \ele_F(\overline{M},N)={}^*F(M)\ot N,
 \end{split}
\end{equation}
 equipped with pre-balancing
 \begin{equation}\label{defin-gamma} \begin{split}  \gamma^X_{M,N}: \ere_F(M,X\triangleright N)\to \ere_F(X^*\triangleright M,N),  \\
 \gamma^X_{M,N} =((\phi^{l}_{X^*,M})^{-1}\ot \id_{F(N)})  (\id_{{}^*M} \ot c_{X,N}),
 \end{split}
\end{equation}
 \begin{equation}\label{defin-eta} \begin{split}  \eta^X_{M,N}: \ele_F(M,X\triangleright N)\to \ele_F(X^*\triangleright M,N),  \\
\eta^X_{M,N} = {}^{*}(c_{X^*,M}) (\phi^{l}_{X^*,F(M)})^{-1}\ot \id_N.
 \end{split}
\end{equation}
Here we are omitting the isomorphisms $ X\simeq {}^*(X^*)$, for any $X\in \ca$, and isomorphisms $\phi^l$ are those presented in \eqref{duals-tensor-product}.

\begin{lema}\label{left-right-SD-modend} Let $A, B\in \ca$ be algebras such that module categories $\ca_A, \ca_B$ are exact. Let $(F, c)\in \Fun_\ca(\ca_A,\ca_B)$ be a $\ca$-module functor. The following statements hold.
\begin{itemize}
\item[(i)] There exists an equivalence of categories ${}^*( -)\boxtimes_\ca \Id: \ca^{\op}_A \botc \ca_B\to {}_A\ca \botc \ca_B$ such that
$$({}^*( -)\boxtimes_\ca \Id)\circ \botc \simeq \botc \circ ({}^*(- )\times \Id)  $$
as $\ca$-balanced functors.

\item[(ii)] If the module end $ \oint_{M\in \ca_A}  (S_F, \beta) $ exists,  then  
$$\widehat{\ot}\circ ({}^*( -)\boxtimes_\ca \Id)\big(\oint_{M\in \ca_A}  (S_F, \beta)\big) \simeq  \oint_{M\in \ca_A}  (\ere_F, \gamma).$$

\item[(iii)]  If the module end $ \oint_{M\in \ca_A}  (D_F, \delta) $ exists,  then  
$$\widehat{\ot}\circ ({}^*( -)\boxtimes_\ca \Id)\big(\oint_{M\in \ca_A}  (D_F, \delta)\big) \simeq  \oint_{M\in \ca_A}  (\ele_F, \eta).$$

\end{itemize}
Here $\widehat{\ot}: {}_A\ca\boxtimes_{\ca} \ca_B \to  {}_A\ca_B$ is the induced  functor from the tensor product, that we have presented in Proposition \ref{basic-stuff} (ii).
\end{lema}
\pf Exactness of module categories $\ca_A, \ca_B$ is needed to ensure existence of functors $S_F, D_F$, see \cite[Thm. 3.3 (4)]{DSS}. Part (i) follows since $\botc \circ ({}^*(- )\times \Id)$  is a $\ca$-balanced functor. There are isomorphisms of $\ca$-balanced functors
\begin{align*} \widehat{\ot}\circ ({}^*( -)\boxtimes_\ca \Id)\circ S_F& \simeq \widehat{\ot}\circ ({}^*( -)\boxtimes_\ca \Id)\circ \botc\circ  (\Id \times F)\\
&\simeq  \widehat{\ot}\circ  \botc\circ  ({}^*( -) \times F)\\
&\simeq \ot \circ  ({}^*( -) \times F)= \ere_F.
\end{align*}
The first isomorphism follows by the definition of $S_F$, the second isomorphism follows from part (i), and the third isomorphism is the one presented in Proposition  \ref{basic-stuff} (ii). Now, part (ii) follows by applying Proposition \ref{equivalences-mod-end} (i). The proof of (iii) follows similarly. 
\epf

\begin{teo}\label{inv-L} Let $A, B\in \ca$ be algebras such that $\ca_A, \ca_B$ are exact module categories. The functor $$\Upsilon: \Fun_\ca(\ca_A, \ca_B)\to \ca^{\op}_A\boxtimes_\ca \ca_B,$$ given by
$$ \Upsilon(F)= \oint_{M\in \ca_A} (S_F,\beta)= \oint_{M\in \ca_A}  \overline{M}\boxtimes_\ca F(M)$$
is well-defined and is a quasi-inverse of the functor $\widetilde{L}$.
\end{teo}
\pf  Recall the definition of the functor $R$ given in Proposition \ref{basic-stuff}, $R:   {}_A\ca_B\to \Fun_\ca(\ca_A,\ca_B),$  $R(V)(X)=X\ot_A V$. It follows from Lemma \ref{aboutL} that, the composition of functors
$$ (\ca_A)^{\op}\boxtimes_{\ca} \ca_B \xrightarrow{^*( -) \boxtimes\Id} {}_A\ca\boxtimes_{\ca} \ca_B \xrightarrow{\ot}  {}_A\ca_B \xrightarrow{ R } \Fun_\ca(\ca_A,\ca_B)$$
is isomorphic to $\widetilde{L}$. Thus, it is enough to show that, the functor 
$$ \Psi: \Fun_\ca(\ca_A,\ca_B)\to {}_A\ca_B,$$
given by 
\begin{equation}\label{quasi-invR-psi}
\Psi(F)= \oint_{M\in \ca_A}  (\ere_F, \gamma)= \oint_{M\in \ca_A}  {}^*M\otimes F(M)
\end{equation}
is well-defined and it is a quasi-inverse of $R$. Since we know that $R$ is an equivalence, we denote by $\Psi$ an adjoint equivalence to $R$. Take $F\in \Fun_\ca(\ca_A,\ca_B)$, and $V\in {}_A\ca_B$, then
\begin{align*} \Hom_{(A,B)}(V, \Psi(F)) &\simeq \Nat_m(R(V), F)\\
&\simeq  \oint_{M\in \ca_A}  (\Hom_{B} (M\ot_A V, F(M)), \beta)\\
& \simeq \oint_{M\in \ca_A} ( \Hom_{ (A,B)} ( V,  {}^*M\otimes F(M) ),\delta)\\
& \simeq  \Hom_{ (A,B)} ( V,  \oint_{M\in \ca_A}  {}^*M\otimes F(M))
\end{align*}
The second isomorphism follows from Proposition \ref{end-natural-transf}. Here, the isomorphism $\beta$ is the one described in \eqref{beta-for-homs}. The third isomorphism follows from Lemma \ref{dual-module1} (ii); one can easily verify that if
$$\delta^X_{M,N}:  \Hom_{ (A,B)} ( V,  {}^*M\otimes F(X\ot N))\to  \Hom_{ (A,B)} ( V,  {}^*M\otimes X\ot F(N))$$
 is defined as $\delta^X_{M,N}(h)=(\id_{{}^*M}\ot c_{X,N})\circ h$, then the naturality of $\Phi$ implies that
$$ \delta^X_{M,N}(\Phi(\alpha))=\Phi(\beta^X_{M,N}(\alpha)), $$
for any $\alpha\in \Hom_{B} (M\ot_A V, F(X\ot N))$. Here 
$$\Phi:  \Hom_{B} (M\ot_A V, F(X\ot N))\to  \Hom_{ (A,B)} ( V,  {}^*M\otimes F(X\ot N))$$
is the natural isomorphism described in \eqref{def-phi-b}. Thus, the third isomorphism follows by applying Proposition \ref{equivalences-mod-end} (i).
The last isomorphism follows from Proposition \ref{equivalences-mod-end} (ii). 
\epf

As an immediate consequence of the above Theorem, we have the following results.
\begin{cor}\label{pw-bim} Let $A\in\ca$ be an algebra such that $\ca_A$ is an exact module category. There is an isomorphism of $A$-bimodules
$$ A\simeq \oint_{M\in \ca_A}  {}^*M\otimes M.$$\qed
\end{cor}

\begin{cor}\label{psi-phi-inv-eachother}  Let $\Mo$, $\No$ be exact indecomposable $\ca$-module categories. 
 If $U\in \Mo$, $V\in \No$ and $F\in \Fun_\ca(\Mo, \No)$, there are isomorphisms
\begin{equation}\label{pp1}
 \oint_{M\in \Mo}  \widetilde{L}_{\Mo,\No}(\overline{M}\boxtimes_\ca F(M)) \simeq F,
\end{equation} 
\begin{equation}\label{pp2}
 \oint_{M\in \Mo}  \overline{M}  \boxtimes_\ca \widetilde{L}_{\Mo,\No}(\overline{U}\boxtimes_\ca V)(M) \simeq  \overline{U} \boxtimes_\ca V. \end{equation} \qed
\end{cor}

\begin{rmk} If $\ca=\vect_\ku$, Corollary \ref{pw-bim} reduces to \cite[Corollary 2.9]{FSS0}.
\end{rmk} 

In \cite[Lemma 3.8]{FSS0} it was proven that, for a right exact functor $F:\Mo\to \No$, where $\Mo, \No$ are abelian categories, there is an isomorphism 
$$\int_{N\in \No} \overline{F^{r.a.}(N)}\boxtimes  N \simeq   \int_{M\in \Mo} \overline{M} \boxtimes F(M). $$
 The next result is a generalization of that result; essentially it says that, for a $\ca$-module functor $F:\Mo\to \No$, there is an isomorphism 
$$\oint_{N\in \No} \overline{F^{r.a.}(N)}\btc N \simeq   \oint_{M\in \Mo} \overline{M}\btc F(M).$$
The proof, however,  is more complicated than the proof of \cite[Lemma 3.8]{FSS0}, since in module ends there is a new ingredient (the pre-balancing $\beta$) that has to be taken into account.
\begin{prop} Let $\Mo, \No$ be exact indecomposable left $\ca$-module categories. Assume that, $(F,c)\in \Fun_\ca(\Mo, \No)$ is a module functor with right adjoint $(F^{r.a.},d)\in \Fun_\ca(\No, \Mo)$. Recall the functors $D_F, S_F$ defined in \eqref{defin-df}, \eqref{defin-sf} with their pre-balancings $\delta, \beta$. There is an isomorphism
\begin{equation}\label{end-adj} \oint_{N\in \No} (\overline{F^{r.a.}(N)}\btc N, \delta) \simeq   \oint_{M\in \Mo} (\overline{M}\btc F(M),\beta).
\end{equation}
\end{prop}
\pf  Since $\Mo, \No$ are exact indecomposable, we can assume that, there are algebras $A, B\in \ca$ such that $\Mo=\ca_A, \No=\ca_B$. Using Lemma \ref{left-right-SD-modend} (ii), (iii), it will be enough to prove that there are isomorphisms
$$ \oint_{N\in \ca_B}  (\ele_{F^{r.a.}}, \eta)\simeq \oint_{M\in \ca_A}  (\ere_F, \gamma),$$
as objects in $  {}_A\ca_B$. Here $\eta, \gamma$ are defined in \eqref{defin-eta}, \eqref{defin-gamma}. Since the functor $R:   {}_A\ca_B\to \Fun_\ca(\ca_A,\ca_B),$  $R(V)(X)=X\ot_A V$ is an equivalence of categories, using Proposition \ref{equivalences-mod-end} (iii), it will be enough to prove that, there is an isomorphism
\begin{equation}\label{equi-prop11}
 \oint_{N\in \ca_B}  \big(R( {}^*F^{r.a.}(N)\ot N), R(\eta)\big) \simeq \oint_{M\in \ca_A}  \big( R({}^*M \ot F(M)), R(\gamma)\big).
\end{equation}
Since the functor $R$ is a quasi-inverse of the functor $\Psi:\Fun_\ca(\ca_A,\ca_B)\to {}_A\ca_B$, presented in \eqref{quasi-invR-psi}, it follows that 
$$\oint_{M\in \ca_A}  \big( R({}^*M \ot F(M)), R(\gamma)\big)\simeq F,$$
and
$$\oint_{N\in \ca_B}  \big( R({}^*N \ot N), R(\gamma)\big) \simeq \Id_{\ca_B}.$$
Hence, to prove isomorphism \eqref{equi-prop11} of functors, it is sufficient to prove that, there is an isomorphism
$$  \oint_{N\in \ca_B}  \big(R( {}^*F^{r.a.}(N)\ot N), R(\eta)\big)(U)\simeq  \oint_{M\in \ca_B}  \big( R({}^*M \ot M), R(\gamma)\big) (F(U))$$
For any $U\in \ca_A$. Applying Theorem \ref{parameter-end}, it will be enough to prove that, there is an isomorphism
$$ \oint_{N\in \ca_B} \big(U\ot_A {}^*F^{r.a.}(N)\ot N , \widehat{\eta}\big) \simeq \oint_{M\in \ca_B}  \big( F(U)\ot_B {}^*M \ot M), \widehat{\gamma}\big),$$
where 
$$\widehat{\eta}^X_{M,N}= R(\eta)_U=\id_U\ot_A  {}^*(d_{X^*,M})(\phi^l_{X^*,F^{r.a.}(M)})^{-1}\ot \id_N, $$
$$ \widehat{\gamma}^X_{M,N}=R(\gamma)_{F(U)}=\id_{F(U)}\ot_B (\phi^l_{X^*,M})^{-1}\ot \id_N,$$
for any $X\in \ca,$ $M, N\in \ca_B$. For this purpose, we shall construct natural isomorphisms $$a_{U,M}:F(U)\ot_B {}^*M \to U\ot_A {}^*F^{r.a.}(M) $$ such that 
\begin{equation}\label{iso-end-a1} 
\widehat{\eta}^X_{M,N}(a_{U,M}\ot \id_{X\ot N})= (a_{U, X^*\ot M}\ot \id_N)\widehat{\gamma}^X_{M,N}.
\end{equation}
It will follow then from Proposition \ref{equivalences-mod-end} (i) the desired isomorphism between module ends, and this will finish the proof of the Proposition.
\medbreak 

Recall the isomorphisms $\Phi^A_{M,X,N}, \Psi^A_{M,X,N}$ defined in \eqref{defini-isos2}, \eqref{defini-isos22}
\begin{equation*}\begin{split}
 \Phi^A_{M,X,N}: \Hom_A(M, X\ot N)\to \Hom_\ca(M\ot_A {}^*N, X),\\
  \Phi^A_{M,X,N}(\alpha)\pi^A_{M, {}^*N}=(\id_X\ot \ev_N)(\alpha\ot\id_{ {}^*N}),
 \end{split}
\end{equation*}
\begin{equation*}\begin{split}
  \Psi^A_{M,X,N}:\Hom_\ca(M\ot_A {}^*N, X)\to  \Hom_A(M, X\ot N),\\
  \Psi^A_{M,X,N}(\alpha)=(\alpha \pi^A_{M, {}^*N}\ot\id_N)(\id_M\ot \coev_N).
    \end{split}
\end{equation*}
We shall also denote natural isomorphisms 
$$\omega_{M,N}: \Hom_B(F(M),N)\to \Hom_A(M, F^{r.a.}(N)),  $$
comming from the adjunction $(F,  F^{r.a.})$. Naturality of $\omega$ implies that for any morphism $f:N\to \widetilde{N}$ in $\ca_B$, and any $\alpha\in \Hom_B(F(M),N)$ we have that 
\begin{equation}\label{naturality-omega} \omega_{M,\widetilde{N}}(f \, \alpha)=F^{r.a.}(f)  \omega_{M,N}(\alpha).
\end{equation}
This equation implies in particular that
\begin{equation}\label{naturality-omega-part} 
 \omega_{U, Y\ot N}( \Psi^B_{F(U),Y,M}(\id))=F^{r.a.}( \Psi^B_{F(U),Y,M}(\id))  \omega_{U,F(U)}(\id).
\end{equation}
Define isomorphisms 
$$a_{U,M}:F(U)\ot_B {}^*M \to U\ot_A {}^*F^{r.a.}(M) $$ 
induced by the natural isomorphisms
$$\Hom_\ca(F(U)\ot_B {}^*M, Z )\xrightarrow{  \Psi^B} \Hom_B(F(U),Z\ot M) \xrightarrow{  \omega} \Hom_A(U,F^{r.a.}(Z\ot M))$$ $$ \xrightarrow{  } \Hom_A(U, Z\ot F^{r.a.}( M)) \xrightarrow{  \Phi^A} \Hom_\ca(U\ot_A  {}^*F^{r.a.}( M), Z), $$
for any $Z\in \ca$. This means that
$$a^{-1}_{U,M}=\Phi^A_{U,Y, F^{r.a.}(M)}\big( d_{Y,M}\, \omega_{U,Y\ot N}( \Psi^B_{F(U),Y,M}(\id) )\big), $$ where $Y=F(U)\ot_B {}^* M.$
Using the definition of $\Phi^A$ one gets that
\begin{align}\label{about-aphi}\begin{split}
a^{-1}_{U,M} \pi^A_{U, {}^*F^{r.a.}(M)}&=(\id_Y\ot \ev_{F^{r.a.}(M)})\\& \big( d_{Y,M}\, \omega_{U,Y\ot M}( \Psi^B_{F(U),Y,M}(\id)) \ot \id_{ {}^*F^{r.a.}(M)}\big).
\end{split}
\end{align}
Here we are again denoting $Y=F(U)\ot_B{}^*M$. Equation \eqref{iso-end-a1} is equivalent to 
\begin{align*} \begin{split}
(a^{-1}_{U, X^*\ot M}\ot \id_N) \widehat{\eta}^X_{M,N} &(\pi^A_{U, {}^*F^{r.a.}(M)}\ot \id_{X\ot N})=\\ &= \widehat{\gamma}^X_{M,N}(a^{-1}_{U,M}\ot \id_{X\ot N})(\pi^A_{U, {}^*F^{r.a.}(M)}\ot \id_{X\ot N}),
\end{split}
\end{align*}
which in turn (forgetting the last $\id_N$) is equivalent to 
\begin{align}\label{iso-end-a2} \begin{split}
(a^{-1}_{U, X^*\ot M}& \pi^A_{U, {}^*F^{r.a.}(X^*\ot M)} )(\id_U\ot  {}^*(d_{X^*,M})(\phi^l_{X^*,F^{r.a.}(M)})^{-1})=\\ &=(\id_{F(U)}\ot_B (\phi^l_{X^*,M})^{-1})(a^{-1}_{U,M} \pi^A_{U, {}^*F^{r.a.}(M)}\ot \id_{X}).
\end{split}
\end{align}

Using \eqref{about-aphi},  the right hand side of \eqref{iso-end-a2} is equal to
\begin{align*}
&=(\id_{F(U)}\ot_B (\phi^l_{X^*,M})^{-1})(a^{-1}_{U,M}\pi^A_{U, {}^*F^{r.a.}(M)}\ot \id_{X})\\
&= (\id_{F(U)}\ot_B (\phi^l_{X^*,M})^{-1}) \big( \id_{F(U)\ot_B {}^*M}\ot \ev_{F^{r.a.}(M)}\ot\id_{X} \big)   \big(  d_{F(U)\ot_B {}^*M, M}\\
&  \omega_{U, F(U)\ot_B {}^*M\ot M} (\Psi^B_{F(U), F(U)\ot_B {}^*M,  M}(\id)) \ot \id_{{}^*F^{r.a.}(M)\ot X }\big)\\
&= (\id_{F(U)}\ot_B (\phi^l_{X^*,M})^{-1}) \big( \id_{F(U)\ot_B {}^*M}\ot \ev_{F^{r.a.}(M)}\ot\id_{X} \big)   \big(  d_{F(U)\ot_B {}^*M, M}\\
& F^{r.a.}(\Psi^B_{F(U), F(U)\ot_B {}^*( M),  M}(\id)) \,  \omega_{U, F(U)} (\id)\ot \id_{{}^*F^{r.a.}(M)\ot X }\big)
\end{align*}
The last equality follows from \eqref{naturality-omega-part}. It follows from \eqref{about-aphi}, that the left hand side of \eqref{iso-end-a2} is equal to
\begin{align*} &=\big(\id_{F(U)\ot_B {}^*(X^*\ot M)}\ot \ev_{F^{r.a.}(X^*\ot M)} \big) \big(d_{F(U)\ot_B {}^*(X^*\ot M), X^*\ot M} \\
& \omega_{U, F(U)\ot_B {}^*(X^*\ot M)\ot X^*\ot M} (\Psi^B_{F(U), F(U)\ot_B {}^*(X^*\ot M), X^*\ot M}(\id)) \ot \id_{{}^*F^{r.a.}(X^*\ot M)}\big)\\& \big( \id_U\ot  {}^*(d_{X^*,M}) 
(\phi^l_{X^*,F^{r.a.}(M)})^{-1}\big) \\
&=\big(\id_{F(U)\ot_B {}^*(X^*\ot M)}\ot \ev_{F^{r.a.}(X^*\ot M)} \big)\big( \id_{F(U) \ot_B {}^*(X^* \ot M) \ot F^{r.a}(X^* \ot M)} \ot \\ &\ot  {}^*(d_{X^*,M}) 
(\phi^l_{X^*,F^{r.a.}(M)})^{-1}\big)\\ &\big(d_{F(U)\ot_B {}^*(X^*\ot M), X^*\ot M}\, F^{r.a.}(\Psi^B_{F(U), F(U)\ot_B {}^*(X^*\ot M),X^*\ot M}(\id) )\ot \id \big)\\
&\big(\omega_{U, F(U)} (\id) \ot \id\big)\\
&= \big(\id_{F(U)\ot_B {}^*(X^*\ot M)}\ot \ev_{X^*\ot F^{r.a.}( M)} (d_{X^*,M}\ot (\phi^l_{X^*,F^{r.a.}(M)})^{-1} ) \big)\\ & \big(d_{F(U)\ot_B {}^*(X^*\ot M), X^*\ot M}\, F^{r.a.}(\Psi^B_{F(U), F(U)\ot_B {}^*(X^*\ot M),X^*\ot M}(\id) \ot \id \big)\\
&\big(\omega_{U, F(U)} (\id) \ot \id\big)\\
&= \big(\id_{F(U)\ot_B {}^*(X^*\ot M)}\ot \ev_{X^*\ot F^{r.a.}( M)}  (\id \ot (\phi^l_{X^*,F^{r.a.}(M)})^{-1})\big)\\& \big( d_{F(U)\ot_B {}^*(X^*\ot M) \ot X^*, M} \ot \id_{{}^*F^{r.a.}(M)\ot X} \big) \\  &\big(  F^{r.a.}(\Psi^B_{F(U), F(U)\ot_B {}^*(X^*\ot M),X^*\ot M}(\id) \ot \id \big)
\big(\omega_{U, F(U)} (\id) \ot \id\big).
\end{align*}
The second equation follows from \eqref{naturality-omega-part}, the third equality follows from \eqref{duality-morph}, the fourth equality follows from \eqref{modfunctor1} for the module functor $(F^{r.a.},d) $, which in this case implies that 
$$(\id \ot d_{X^*,M}) d_{F(U)\ot_B {}^*(X^{*}\ot M), X^{*}\ot M}=d_{F(U)\ot_B {}^*(X^{*}\ot M)\ot X^*, M}.$$
At last, using the definition of $\Psi^B$, \eqref{ev-tensor-prod} and the rigidity axioms one can see that 
\begin{align*}
 &\big(\id_{F(U)\ot_B {}^*(X^*\ot M)}\ot \ev_{X^*\ot F^{r.a.}( M)}  (\id\ot (\phi^l_{X^*,F^{r.a.}(M)})^{-1})\big)\\ &\big( d_{F(U)\ot_B {}^*(X^*\ot M) \ot X^*, M} \ot \id\big)   \big(  F^{r.a.}(\Psi^B_{F(U), F(U)\ot_B {}^*(X^*\ot M),X^*\ot M}(\id) ) \ot \id \big)\\& (\omega_{U,F(U)}(\id) \ot \id_{{}^*F^{r.a}(M) \ot X}) \\
 &= (\id_{F(U)}\ot_B (\phi^l_{X^*,M})^{-1}) \big( \id_{F(U)\ot_B {}^*M}\ot \ev_{F^{r.a.}(M)}\ot\id_{X} \big) \big(  d_{F(U)\ot_B {}^*M, M}\\
& F^{r.a.}(\Psi^B_{F(U), F(U)\ot_B {}^*( M),  M}(\id))(\omega_{U,F(U)}(\id) \ot \id_{{}^*F^{r.a}(M) \ot X}).
\end{align*}
This implies \eqref{iso-end-a2}, and finishes the proof of the Proposition.
\epf

\subsection{A formula for the relative Serre functor}

Let $\Mo$, $\No$ be  exact indecomposable left $\ca$-module categories, and recall the functors
$$L=L_{\Mo,\No}: \Mo^{\op}\boxtimes_\ca \No \to \Fun_\ca(\Mo^{\bop}, \No),$$
$$ \widetilde{L}=\widetilde{L}_{\No,\Mo^{\bop}}: \No^{\op}\boxtimes_\ca \Mo^{\bop} \to \Fun_\ca(\No, \Mo^{\bop})$$
described in \eqref{L-equivalence1} and \eqref{L-equivalence2}. Note that subindices of $ \widetilde{L}$ are different to those presented in \eqref{L-equivalence2}.

\begin{lema}\label{adjoint-L} Use the above notation. For any $M\in \Mo$, $N\in \No$ there exists an equivalence of module functors
\begin{equation}\label{left-adj-L} L_{\Mo,\No}(\overline{M}\btc N)^{l.a.} \simeq \widetilde{L}_{\No,\Mo^{\bop}}(\overline{N}\btc  \overline{\overline{M}})
\end{equation}
\end{lema}
\pf If $\Bc$ is an exact indecomposable right $\ca$-module category, define the auxiliary functors
$$H^{\Bc}_B:\Bc^{\op}\to \ca,\quad  \ere^\No_N:\ca\to \No,$$
$$H^{\Bc}_B= \uhom_\Bc( -, B), \quad \ere^\No_N= - \triangleright N,$$
for any $B\in \Bc$, $N\in \No$. A straightforward computation shows that
$$(H^{\Bc}_B)^{l.a.}(X)=B\triangleleft X^*, \quad (\ere^\No_N)^{l.a.}(N')={}^*\uhom_\No( N', N) $$
for any $X\in \ca$, $N'\in \No$. Since $ L(\overline{M}\btc N)=\ere^\No_N \circ H^{\Mo^{\op}}_{\overline{M}}$, then
\begin{align*}
L(\overline{M}\btc N)^{l.a.}&\simeq  (H^{\Mo^{\op}}_{\overline{M}})^{l.a.}\circ (\ere^\No_N)^{l.a.}\\
&\simeq\overline{M}\triangleleft \uhom_\No( -, N)=\uhom_\No( -, N)^* \triangleright M\\
&\stackrel{\eqref{hom-op}}\simeq \uhom_{\No^{\op}}( -, \overline{N})^{***} \triangleright M=\widetilde{L}(\overline{N}\btc  \overline{\overline{M}}).
\end{align*} 
In the second equivalence, we are using the canonical isomorphisms ${}^*X^*\simeq X$. 
\epf

The next result is a formula for the relative Serre functor similar to the formula for the Nakayama functor given in \cite{FSS}. Let $\Mo$ be an exact indecomposable left $\ca$-module category. Let us denote by $\blacktriangleleft :\Mo^{\op}\times \ca\to \Mo^{\op}$ the action of the opposite module category, that is, the one determined by 
\begin{equation}\label{aop}
\overline{M}\blacktriangleleft X=\overline{X^*\triangleright M},
\end{equation}
 for any $M\in \Mo$, $X\in\ca$.  For any $M\in \Mo$ the functor 
$$T_M:(\Mo^{\op})^{\op}\times \Mo^{\op}\to \Mo,$$
$$T_M(U,V)=\uhom_\Mo(M,V)^*\triangleright U, $$
 has a pre-balancing 
$$ \gamma^X_{U,V}: T_M(U\blacktriangleleft X, V)\to T_M(U, V\blacktriangleleft {}^* X),$$
given as the composition
$$T_M(U\blacktriangleleft X, V)\! =\! \uhom_\Mo(M,V)^*\triangleright (X^*\triangleright U)\xrightarrow{m^{-1}} ( \uhom_\Mo(M,V)^* \ot X^*)\triangleright U$$
$$\xrightarrow{ } (X\ot\uhom_\Mo(M,V))^* \triangleright U \to \uhom_\Mo(M,X \triangleright V)^* \triangleright U =T_M(U, V\blacktriangleleft {}^* X).$$
Thus we can consider the coend
$$\oint^{\overline{U}\in \Mo^{\op}} (T_M,\gamma). $$
Since $T$ can be thought of as a functor $T:(\Mo^{\op})^{\op}\times \Mo^{\op}\to \Fun(\Mo,\Mo)$, $T(U,V)(M)= T_M(U,V)$, then using the parameter theorem described in Section \ref{parameter-coends}, we have a functor
$$M\mapsto \oint^{\overline{U}\in \Mo^{\op}} (T_M,\gamma).$$
We shall denote this functor as 
$$  \oint^{\overline{U}\in \Mo^{\op}} (T_{-},\gamma)=\oint^{\overline{U}\in \Mo^{\op}} (\uhom_{\Mo}(- , U)^*\triangleright U,\gamma).$$
It follows from Lemma \ref{internal-hom-actions} that, $ \oint^{\overline{U}\in \Mo^{\op}} (T_{-},\gamma):\Mo\to \Mo^{\bop}$ is a $\ca$-module functor.
\begin{teo}\label{formula-serre-relat} Let $\Mo$ be an exact indecomposable left $\ca$-module category.  There exists an equivalence of $\ca$-module functors
\begin{equation}\label{formula-serre}
\sy_{\Mo} \simeq \oint^{\overline{U}\in \Mo^{\op}} (\uhom_{\Mo}(- , U)^*\triangleright U,\gamma),
\end{equation}  
\end{teo}
 \pf Let $\Mo,\No$ be a pair of exact indecomposable left $\ca$-module categories.  To prove the expression for the relative Serre functor, we will first compute a quasi-inverse of the functor $L=L_{\Mo,\No}$ and then use equivalence \eqref{eles-serre}. 
 \medbreak
 
 Recall that $\blacktriangleleft :\Mo^{\op}\times \ca\to \Mo^{\op}$ is the action of the opposite module category. That is,
 $$ \overline{M}\blacktriangleleft X=\overline{X^*\triangleright M},$$
 for any $M\in \Mo$, $X\in\ca$.
Let us denote by $D:(\Mo^{\op}\btc \No)^{\op}\to \No^{\op}\btc \Mo^{\bop}$, the functor determined by $D(\overline{\overline{M}\btc N})=\overline{N}\btc \overline{\overline{M}}$, $M\in \Mo$, $N\in \No.$ The functor $D$ is an equivalence of categories.

Take $(F, c)\in \Fun_\ca(\Mo^{\bop},\No)$. This means that, we have isomorphisms $c_{X,M}: F(X^{**} \triangleright M) \to X \triangleright F(M)$, for any $M\in \Mo, X\in \ca$. Define
$$\Theta_F: (\Mo^{\op})^{\op}\times \Mo^{\op}\to \Mo^{\op}\btc \No,$$
$$ \Theta_F(\overline{\overline{U}},\overline{V})=\overline{V}\btc F(\overline{\overline{U}}).  $$
The functor $\Theta_F$ has a pre-balancing
$$\nu^X_{U,V}: \Theta_F(\overline{\overline{U}} \blacktriangleleft X, \overline{V}) \to \Theta_F(\overline{\overline{U}}, \overline{V} \blacktriangleleft {}^*X),$$
$$\nu^X_{U,V}=b^{-1}_{V, {}^*X, F(U)} (\id_V\btc c_{{}^*X,U}). $$
Here $b_{V,X,U}: V \blacktriangleleft X \btc U \to V\btc X  \triangleright  U$ is the balancing of the $\ca$-balanced functor $\btc$. Define $\chi:\Fun_\ca(\Mo^{\bop},\No)\to \Mo^{\op}\btc \No$ the functor given by 
$$\chi(F)=\oint^{\overline{U}\in \Mo^{\op}} (\Theta_F, \nu)=\oint^{\overline{U}\in \Mo^{\op}} \overline{U}\btc F(\overline{\overline{U}}). $$
The existence of these coends follows  from the existence of the ends presented in Theorem \ref{inv-L} and the relation between ends and coends for left and right module categories given in Lemma \ref{left-right}.  Let us prove that $\chi$ is a quasi-inverse of $L$. Since we already know that $L$ is a category equivalence, it is enough to prove that
$$\chi(L(\overline{M}\btc N)) \simeq \overline{M}\btc N$$
for any $M\in \Mo, N\in \No$. Since $D$ is a category equivalence, this is equivalent to prove that 
\begin{equation}\label{transpo-delig}
D(\overline{ \chi(L(\overline{M}\btc N))}) \simeq D(\overline{ \overline{M}\btc N})=\overline{N}\btc \overline{\overline{M}}.
\end{equation} 
for any $M\in \Mo, N\in \No$. The left hand side of \eqref{transpo-delig} is equal to
\begin{align*}
D(\overline{ \chi(L(\overline{M}\btc N))})&=D\big(\overline{\oint^{\overline{U}\in \Mo^{\op}} \overline{U}\btc L(\overline{M}\btc N)(\overline{\overline{U}})}\big)\\
&\simeq D\big( \oint_{\overline{\overline{U}}\in \Mo^{\bop}} \overline{\overline{U}\btc L(\overline{M}\btc N)(\overline{\overline{U}})}\big)\\
&\simeq \oint_{\overline{\overline{U}}\in \Mo^{\bop}} \overline{ L(\overline{M}\btc N)(\overline{\overline{U}})} \btc \overline{\overline{U}}\\
&\simeq\oint_{V\in \No}   \overline{V}\btc L(\overline{M}\btc N)^{l.a.}(V)\\
&\stackrel{\eqref {left-adj-L}}\simeq\oint_{V\in \No}   \overline{V}\btc \widetilde{L}(\overline{N}\btc  \overline{\overline{M}})(V)\stackrel{\eqref{pp2}}\simeq \overline{N}\btc  \overline{\overline{M}}.
\end{align*}
The first isomorphism follows from Lemma \ref{right-left}, the second one follows from Proposition \ref{equivalences-mod-end} (iii), and the third isomorphism follows from Proposition \ref{end-adj}.

\smallbreak

Taking $\No=\Mo^{\bop}$ and using \eqref{eles-serre}, it follows that
$$\widetilde{L}_{\Mo,\Mo^{\bop}}(\chi(\Id))\simeq L_{\Mo,\Mo^{\bop}}(\chi(\Id))\circ 	\sy_\Mo\simeq \sy_\Mo,$$
and  we obtain the desired description of the relative Serre functor.
 \epf

\begin{rmk} If $\ca=\vect_\ku$ and $\Mo$ is a semisimple category, the (right) relatvie Serre functor coincides with the (right) Nakayama functor. In this case, formula \eqref{formula-serre} coincides with the formula for the (right) Nakayama functor presented in \cite[Definition 3.14]{FSS0}.
\end{rmk}

\subsection{Correspondence of module categories for Morita equivalent tensor categories}

Assume that $\ca, \Do$ are  Morita equivalent tensor categories. This means that, there is an invertible exact $(\ca, \Do)$-bimodule category $\Bc$. We can assume that $\Do=\End_\ca(\Bc)^{\rev},$ and the right action of $\Do$ on $\Bc$ is given by evaluation
$$ \triangleleft:\Bc\times \Do\to \Bc,\,\,  B\triangleleft F=F(B).$$
It was proven in \cite[Theorem 3.31]{EO} that, the maps
$$ \Mo \mapsto \Fun_\ca(\Bc, \Mo), \quad \No \mapsto \Fun_{\Do}(\Bc^{\op}, \No)$$
are bijections, one the inverse of the other, between equivalence classes of exact $\ca$-module categories and exact $\Do$-module categories. We shall give another proof of this fact by showing  an explicit equivalence of $\Do$-module categories $$ \No\simeq \Fun_{\ca}(\Bc, \Fun_{\Do}(\Bc^{\op}, \No)),$$ 
for any exact indecomposable $\Do$-module category $\No$.
\medbreak

For any $(H,d)\in \Fun_{\ca}(\Bc, \Fun_{\Do}(\Bc^{\op}, \No)),$ define 
$$S_H: \Bc^{\op}\times \Bc\to \No,\quad S_H(\overline{B}, C)=H(C)(\overline{B}).$$ This functor comes with isomorphisms 
$$ \beta^X_{B,C}: S_H(B, X\triangleright C) \to S_H(X^*\triangleright B, C),$$ 
$$  \beta^X_{B,C} =\big(d_{X,C}\big)_B,$$
for any $X\in \ca$, $B,C\in \Bc$. 
\begin{lema}\label{S_H-balanced} The functor $S_H$ is  a $\ca$-balanced  functor with balancing given by $b_{B,X,C}: S_H(X^*\triangleright B, C)\to S_H(B, X\triangleright C)$, $b_{B,X,C}=\big(d_{X,C}\big)^{-1}_B$. In particular, there exists a right exact functor $\widehat{S}_H: \Bc^{\op}\btc \Bc \to \No$ such that $\widehat{S}_H\circ \btc \simeq S_H$ as $\ca$-balanced functors.
\end{lema}
\pf Since $(H,d)$ is a module functor, the natural isomorphism $d$ satisfy \eqref{modfunctor1}. This axiom implies that $b$ satisfy \eqref{c-balanced}.
\epf 

We can consider the functor 
$$ \Psi: \Fun_{\ca}(\Bc, \Fun_{\Do}(\Bc^{\op}, \No)) \to \No,$$
$$ \Psi(H)= \oint_{B\in \Bc}   (S_H,\beta)= \oint_{B\in \Bc} H(B)(\overline{B}).$$
\begin{prop}\label{psi-correspond-mod} The functor $\Psi$ is well-defined.
\end{prop}
\pf
The existence of the module end $\Psi(H)$ follows from applying the functor $\widehat{S}_H$ to the module end $\oint_{B\in \Bc} \overline{B}\btc B$, whose existence follow from Proposition \ref{inv-L}, and using Proposition \ref{equivalences-mod-end} (iii). 
\epf

\begin{teo}\label{bij-corresp}   Let $\ca$ be a tensor category and $\Bc$ an indecomposable exact left $\ca$-module category. Consider the finite tensor category $\Do=\End_\ca(\Bc)^{\rev}$, and the functor $\widetilde{L}=\widetilde{L}_{\Bc,\Bc}: \Bc^{\op} \btc \Bc\to \End_\ca(\Bc)$ introduced in Section \ref{Subsection:modulefunctors}. Let $\No$ be an exact indecomposable left $\Do$-module category. Define $$\Phi: \No\to \Fun_{\ca}(\Bc, \Fun_{\Do}(\Bc^{\op}, \No)),$$
$$\Phi(N)(B)(\overline{C})=\widetilde{L}(\overline{C}\btc B)\triangleright N, $$
for any $B, C\in \Bc$, $N\in \No$.
The functors $\Phi$ and $\Psi$ are well-defined, and they establish an adjoint equivalence of left $\Do$-module categories 
$$ \No\simeq \Fun_{\ca}(\Bc, \Fun_{\Do}(\Bc^{\op}, \No)).$$
\end{teo}
\pf Take $N\in \No, B\in \Bc$. It follows immediately that, $\Phi(N)$ is a $\ca$-module functor. That $\Phi$ and $\Phi(N)(B)$ are $\Do$-module functors follow from the bimodule structure of the functor $\widetilde{L}$ \eqref{L-is-bimodule1}, \eqref{L-is-bimodule2}. Let us show that the pair of functors $\Phi$, $\Psi$ is an adjoint equivalence. Take $H\in  \Fun_{\ca}(\Bc, \Fun_{\Do}(\Bc^{\op}, \No)),$ $C_1, C_2\in \Bc$, then
\begin{align*}
\Phi(\Psi(H))(C_1)(\overline{C_2})&=\widetilde{L}(\overline{C_2}\btc C_1)\triangleright   \oint_{B\in \Bc} H(B)(\overline{B})\\
&\simeq  \oint_{B\in \Bc} H(B)(\widetilde{L}(\overline{C_2}\btc C_1)^*( B)) \\
&\simeq  \oint_{B\in \Bc} \widehat{S}_H \big( \widetilde{L}(\overline{C_2}\btc C_1)^*( B)  \btc B \big)\\
& \simeq \widehat{S}_H \big( \oint_{B\in \Bc}  \widetilde{L}(\overline{C_2}\btc C_1)^*( B)  \btc B \big)\\
& \simeq \widehat{S}_H \big( \oint_{B\in \Bc}  \overline{B}  \btc   \widetilde{L}(\overline{C_2}\btc C_1)( B) \big)\\
& \simeq \widehat{S}_H (\overline{C_2}\btc C_1)\simeq H(C_1)(\overline{C_2}).
\end{align*}
The first isomorphism follows since $H(B)$ is a $\Do$-module functor, the second isomorphism follows from the definition of $ \widehat{S}_H $ given in Lemma \ref{S_H-balanced}, the third one follows from Proposition \ref{equivalences-mod-end} (iii), the fourth isomorphism follows  from  \eqref{end-adj}, and the fifth isomorphism follows from \eqref{pp2}. 

Now, let us take $N\in \No$, then
\begin{align*}
\Psi(\Phi(N))&= \oint_{B\in \Bc}   \Phi(N)(B)(\overline{B})= \oint_{B\in \Bc} \widetilde{L}( \overline{B}  \btc  B)\triangleright N\\
&\simeq \Id \triangleright N \simeq N.
\end{align*}
The isomorphism follows from \eqref{pp1}. One can verify, in the above proof of $\Phi(\Psi(H))\simeq H$ and in the proof of $\Psi(\Phi(N))\simeq N$, that  each pre-balancing is used properly.
\epf 

\subsection{The double dual tensor category} Let $\Mo$ be an exact indecomposable left $\ca$-module category. Then the \textit{dual tensor category} $\ca^*_\Mo=\End_\ca(\Mo)$ is again a finite tensor category \cite{EO}. The category $\ca^*_\Mo$ acts on $\Mo$ by evaluation: 
$$\ca^*_\Mo\times \Mo\to \Mo, $$
$$ (F, M)\mapsto F(M).$$
The category $\Mo$ is exact indecomposable over $\ca^*_\Mo$, see \cite[Lemma 3.25]{EO}. Whence, we can consider the tensor category $(\ca^*_\Mo)^*_\Mo=\End_{\ca^*_\Mo}(\Mo).$ There is a canonical tensor functor
$$can: \ca\to  (\ca^*_\Mo)^*_\Mo,$$
$$can(X)(M)=X\rhd M, $$
for any $X\in \ca$, $M\in \Mo$. One can see that $can(X)$ is a $\ca^*_{\Mo}$-module functor. It was proven in \cite[Theorem 3.27]{EO} that the functor $can$ is an equivalence of categories. We shall give an expression of  a quasi-inverse of this functor.

\medbreak

Take $(G, d)\in (\ca^*_\Mo)^*_\Mo$. This means that we have natural isomorphisms
$$d_{F, M}: G( F(M))\to F(G(M)), $$
for any $F\in \ca^*_\Mo$, $M\in \Mo$.
Let us denote 
$$ \Hc_{(G, d)}: \Mo^{\op}\times \Mo\to \ca, $$
$$\Hc_{(G, d)}(M, N)=\uhom(M, G(N)).$$
The functor  $ \Hc_{(G, d)}$ has a pre-balancing $\gamma$ (seeing $\Mo$ as a left module category over $\ca^*_\Mo)$. For any $F\in \ca^*_\Mo $ define
$$ \gamma^F_{M,N}:  \Hc_{(G, d)}(M, F(N))\to  \Hc_{(G, d)}(F^{l.a.}(M), N), $$
(Recall that $F^*=F^{l.a.}$) as the composition
$$ \uhom(M, G(F(N))) \xrightarrow{ \uhom(\id, d_{F,N})}   \uhom(M, F(G(N))) \to$$ $$\xrightarrow{ \eqref{uhomadj} }  \uhom(F^{l.a.}(M), G(N)).$$
Explicitly, using \eqref{uhomadj}, this isomorphism is 
$$\gamma^F_{M,N}= \psi^Z_{F^{l.a.}(M),G(N)}\big(\Omega_{Z \triangleright M,G(N)}(\phi^Z_{M,F(G(N))}(\id_Z))b^{-1}_{Z,M} \big)\circ \uhom(\id, d_{F, N}) $$
where $Z=\uhom(M, F(G(N)))$, and isomorphism $b$ is the module structure of the functor $ F^{l.a.}$. Recall  the  isomorphisms presented in \eqref{Hom-interno} 
\begin{equation*}\begin{split}\phi^X_{M,N}:\Hom_{\ca}(X,\uhom(M,N))\to \Hom_{\Mo}(X\triangleright M,N), \\
\psi^X_{M,N}:\Hom_{\Mo}(X\triangleright M,N)\to \Hom_{\ca}(X,\uhom(M,N)),
\end{split}
\end{equation*}
associated to the pair of adjoint functors $( -\triangleright M , \uhom(M, -))$.
\begin{teo}\label{ddual-tensor}  Let $\Mo$ be an exact indecomposable left $\ca$-module category. The functor $\Upsilon:(\ca^*_\Mo)^*_\Mo \to \ca$ given by
$$\Upsilon(G)= \oint_{M\in \Mo}  (\uhom(M, G(M)), \gamma)$$
is well-defined. The pair of functors $(\Upsilon, can)$ is an adjoint equivalence of categories.
\end{teo}
\pf 
We shall prove that,  there are natural isomorphisms 
$$\Nat_{\!m}(can(X),G) \simeq \Hom_\ca(X,\Upsilon(G)).$$ 

Let us fix $X\in \ca$ and $(G, d)\in (\ca^*_\Mo)^*_\Mo$. Using Proposition \ref{end-natural-transf} we have that
\begin{align}\label{natu11} \Nat_{\!m}(can(X),G) \simeq \oint_{M\in \Mo} \big( \Hom_\Mo( X\triangleright M, G(M)), \beta\big).
\end{align}
Recall that $\Mo$ is thought of as a module category over $\ca^*_\Mo$. According to \eqref{beta-for-homs}, the pre-balancing $\beta$  is, in this case,
$$\beta^F_{M,N}:\Hom_\Mo(X\triangleright  M, G(F(N)) )\to \Hom_\Mo(X\triangleright F^{l.a.}(M), G(N)), $$
$$\beta^F_{M,N}(\alpha)=(\ev_F)_{G(N)} F^{l.a.}(d_{F,N}\alpha)b^{-1}_{X,M}. $$
Here $\ev_F: F^{l.a.}\circ F\to \Id$ is the evaluation of the adjoint pair $(F^{l.a.}, F)$. 
If we denote by $\Omega_{M,N}: \Hom_\Mo(M, F(N))\to  \Hom_\Mo(F^{l.a.}(M), N)$  natural isomorphisms, then $(\ev_F)_M= \Omega_{F(M),M}(\id_{F(M)}), $ for any $M\in \Mo$.

Using Proposition \ref{equivalences-mod-end} (ii) we can consider the module end  
$$\oint_{M\in \Mo} \big(\Hom_\ca(X, \uhom(M, G(M)), \widehat{\gamma}\big),$$
where, the pre-balancing in this case is
$$\widehat{\gamma}^F_{M,N}:\Hom_\ca(X, \uhom(M, G(F(N)))\to \Hom_\ca(X, \uhom(F^{l.a.}(M), G(M))),  $$
$$\widehat{\gamma}^F_{M,N}(\alpha)=\gamma^F_{M,N}\circ \alpha.$$
\begin{claim} Isomorphisms
$$\psi^X_{M, G(N)}:\Hom_\Mo(X\triangleright  M, G(N))\to \Hom_\ca(X, \uhom(M,G(N)))   $$
commutes with the pre-balancings, that is
\begin{equation}\label{psi-commutes-preb} 
\widehat{\gamma}^F_{M,N}\circ  \psi^X_{M, G(F(N))}= \psi^X_{F^{l.a.}(M), G(N)}\circ \beta^F_{M,N}.
\end{equation}
\end{claim}
As a consequence of this claim, using Proposition \ref{equivalences-mod-end} (i), we get an isomorphism of module ends
\begin{align*}
&\oint_{M\in \Mo} \big( \Hom_\Mo( X\triangleright M, G(M)), \beta\big) \simeq  \oint_{M\in \Mo} \big(\Hom_\ca(X, \uhom(M, G(M)), \widehat{\gamma}\big)\\
&\simeq \Hom_\ca(X, \oint_{M\in \Mo} \big(\uhom(M, G(M)), \gamma\big) = \Hom_\ca(X,\Upsilon(G)).
\end{align*} 
The second isomorphism follows from Proposition \ref{equivalences-mod-end} (ii). Combining this isomorphism with \eqref{natu11}  we get the result.

\medbreak

 It remains to prove the claim. Naturality of $\psi, \phi$ and $b$ implies that 
\begin{equation}\label{naturalm1} \psi^X_{M,N}(\alpha (h\triangleright \id_M))= 
\psi^Y_{M,N}(\alpha)\circ h,
\end{equation}
 \begin{equation}\label{naturalm2} \uhom(\id,f) \psi^X_{M,N}(\alpha)=\psi^X_{M,N'}(f\circ \alpha),
 \end{equation}
 \begin{equation}\label{naturalm3}
 \phi^X_{M,N}(\alpha\circ h)=\phi^Y_{M,N}(\alpha)(h\triangleright \id_M),
\end{equation}
 \begin{equation}\label{naturalm4}
 b^{-1}_{Y,M} (h\triangleright \id_{F^{l.a.}(M)})=F^{l.a.}(h\triangleright \id_M) b^{-1}_{X,M},
\end{equation}
 for any morphism $h:X\to Y$ in $\ca$ and any $f:N\to N'$, $M, N', N\in \Mo$.
 
 Let $\alpha\in \Hom_\Mo(X\triangleright  M, G(F(N)))$, and $Z=\uhom(M, F(G(N)))$, then the left hand side of \eqref{psi-commutes-preb} evaluated in $\alpha$ is equal to
\begin{align*}&\gamma^F_{M,N}\circ  \psi^X_{M, G(F(N))}(\alpha)=\\ 
& =\psi^Z_{F^{l.a.}(M),G(N)}\big(\Omega_{Z \triangleright M,G(N)}(\phi^Z_{M,F(G(N))}(\id_Z))b^{-1}_{Z,M} \big) \uhom(\id, d_{F, N}) \\ & \psi^X_{M, G(F(N))}(\alpha)\\
&=\psi^Z_{F^{l.a.}(M),G(N)}\big(\Omega_{Z \triangleright M,G(N)}(\phi^Z_{M,F(G(N))}(\id_Z))b^{-1}_{Z,M} \big) \psi^X_{M, F(G(N))}(d_{F, N}\alpha)\\
&= \psi^X_{F^{l.a.}(M),G(N)}\big(\Omega_{Z \triangleright M,G(N)}(\phi^Z_{M,F(G(N))}(\id_Z))b^{-1}_{Z,M}\\
& (\psi^X_{M, F(G(N))}(d_{F, N}\alpha)\triangleright \id_{F^{l.a.}(M)}) \big)\\
&= \psi^X_{F^{l.a.}(M),G(N)}\big(\Omega_{Z \triangleright M,G(N)}(\phi^Z_{M,F(G(N))}(\id_Z))\\
& F^{l.a.}(\psi^X_{M, F(G(N))}(d_{F, N}\alpha)\triangleright \id_M) b^{-1}_{X,M}\big)\\
&= \psi^X_{F^{l.a.}(M),G(N)}\big( \Omega_{F(G(N)),G(N)}(\id) F^{l.a.}(h)b^{-1}_{X,M}\big)
\end{align*}  
The second equality follows from \eqref{naturalm2}, the third equality follows from \eqref{naturalm1}, the fourth follows from \eqref{naturalm4}, the fifth equality follows from the naturality of $\Omega$. In the last equality the map $h$ is
$$h=\phi^Z_{M,F(G(N))}(\id_Z) ( \psi^X_{M, F(G(N))}(d_{F, N}\alpha)\triangleright \id_M). $$
The right hand side of \eqref{psi-commutes-preb} evaluated in $\alpha$ is equal to
$$\psi^X_{F^{l.a.}(M),G(N)}\big( \Omega_{F(G(N)),G(N)}(\id) F^{l.a.}(d_{F, N}\alpha )b^{-1}_{X,M}\big). $$
It remains to observe that $h=d_{F, N}\alpha$, which follows from \eqref{naturalm3}.
\epf

\end{document}